\documentclass[12pt]{article}

\usepackage[leqno]{amsmath}
\usepackage{amsthm,amscd}
\usepackage{amsfonts} %
\usepackage{amssymb} %

\makeatletter
\def\@cite#1#2{{#1\if@tempswa , #2\fi}}
\def\@lbibitem[#1]#2{\item\if@filesw 
      { \def\protect##1{\string ##1\space}\immediate
        \write\@auxout{\string\bibcite{#2}{#1}}}\fi\ignorespaces}
\makeatother
\newtheorem{theorem}{Theorem}[section]
\newtheorem{co}{Corollary}[section]
\newtheorem{lm}{Lemma}[section]
\newtheorem{rem}{Remark}[section]

\newtheorem{definition}{Definition}[section]
\newtheorem{assumption}{Assumption}[section]
\def\Proof{\noindent{\sl Proof.}\qquad}
\def\QED{\hfill\hbox{\vrule width 4pt height 6pt depth 1.5pt}\par\bigskip}

\newtheorem{prop}{Proposition}[section]
\def\diag{\mathop{\rm diag}}

\def\kmms{\kern-\mathsurround}

\newcommand{\trans}{{^t}}

\newcommand{\calX}{{\cal X}}
\newcommand{\calY}{{\cal Y}}
\newcommand{\calZ}{{\cal Z}}
\newcommand{\calG}{{\cal G}}
\newcommand{\calN}{{\cal N}}

\newcommand{\calM}{{\cal M}}

\newcommand{\bbR}{{\mathbb R}}
\newcommand{\bbS}{{\mathbb S}}

\title{Star-shaped distributions and their generalizations }
\author{
Hidehiko Kamiya \\ 
{\it Faculty of Economics, 
Okayama University} \\
Akimichi Takemura \\
{\it Graduate School of Information Science and Technology} \\  
{\it University of Tokyo} \\ 
and \\
Satoshi Kuriki\\
{\it The Institute of Statistical Mathematics}}

\date{}

\begin{document}
\maketitle

\begin{abstract}
\noindent
Elliptically contoured distributions can be considered to be the
distributions for which the contours of the density functions are
proportional ellipsoids.  We generalize elliptically contoured
densities to ``star-shaped distributions'' with 
concentric star-shaped contours  and  show that many
results in the former case continue to hold in the more general case.
We develop a general theory in the framework of abstract group
invariance so that the results can be applied to other cases as well,
especially those involving random matrices.

\end{abstract}

\section{Introduction}

Elliptically contoured distribution is a convenient generalization of
the multivariate normal distribution and now there exists substantial
literature on elliptically contoured distribution, e.g., 
Chapter  2 of the third edition of 
\cite{Anderson2003},
\cite{Fang-Anderson}, \cite{Fang-Zhang90}, and \cite{Gupta-Varga}.
Density $f(x)$ of an elliptically contoured distribution in ${\mathbb R}^p$ 
can be written as $f(x) = f_{\calG}(g(x))$  
with $g(x)=(x\trans \Sigma^{-1} x)^{1/2}$, 
where $\Sigma$ is a $p \times p$ positive definite matrix. 
Under this distribution, the ``length'' $g(x)$ and 
the ``direction'' $x/g(x)$ are independent.   
Moreover, by changing $f_{\calG}( \cdot ),$ 
we can construct the elliptically contoured distribution
with an arbitrary distribution of $g(x).$ 
Because the distribution of $x/g(x)$ is common to all the 
elliptically contoured distributions with the same $\Sigma,$ 
distributional results concerning $x/g(x)$ derived under the 
assumption of normality continue to hold for all elliptically 
contoured distributions having the same $\Sigma.$
This property is often referred to as ``null robustness'' and 
has been extensively discussed in the literature 
(e.g. \cite{Kariya-Sinha}).

However, the class of elliptically contoured distributions seems to  
be too narrow. 
It does not include, e.g., a simple distribution with a density in 
$\bbR^2$ 
whose contours are concentric squares. %
Note that elliptically contoured distributions differ 
from the multivariate normal distribution only in the distribution
of the one-dimensional length.   
Therefore, in the framework of elliptically contoured distributions, 
we can not consider non-normality which is exhibited in skewness or 
asymmetry of distributions.

As a matter of fact, some properties of elliptically contoured 
distributions, including the above-mentioned 
independence of length and direction and the null 
robustness, continue to hold beyond the class of elliptically 
contoured distributions if we define the ``length'' properly.
We %
extend the class of elliptically contoured distributions to 
a class of distributions called star-shaped distributions, 
whose densities have arbitrary star-shaped sets as their contours. 

In this paper, 
star-shaped distribution is developed 
in the general framework of group invariance, especially  
global cross sections and the associated orbital decompositions.
Our primary concern is the star-shaped distribution, but 
the general theory can be applied to problems about 
distributions of random matrices, including the case where 
the group action is non-free.  
For group invariance in statistics in general, see
a recent survey by \cite{Helland}. 
Actually, the star-shaped distributions have also been 
considered under the name ``$v$-spherical distributions'' 
by \cite{Fernandez-Osiewalski-Steel} (see also \cite{Ferreira-Steel}), 
but from a less algebraic point of view.

The material in %
the present 
paper is based on two 
earlier drafts of the authors, \cite{Takemura-Kuriki96},
\cite{Kamiya-Takemura}. 
We give a unified presentation of relevant and original results from 
these drafts in view of the current literature.

The organization of this paper is as follows. 
In Sections 2-3, we develop a general theory 
in the framework of group invariance. 
In Section 2 we study orbital decomposition and global cross 
sections. 
Based on the arguments there, we define decomposable distributions and 
investigate the associated distributional problems in Section 3.
The results in those sections are applied to star-shaped distributions 
in Section 4.  
In Section 5, further applications to random matrices
are presented.  
Some technical details %
are given in the Appendix.

\section{Orbital decomposition and global cross sections}
\label{sec:cross-section}

In this section we review some basic notions about group actions  
and investigate some properties of %
global cross sections.
Our approach is based on global cross sections, but 
there is another approach---the one based on proper actions and 
quotient measures, for which the reader is referred to 
the significant papers \cite{Andersson} and 
Andersson, Br\o ns and Jensen (1983).

\subsection{Orbital decomposition}

Let a group $\calG$ act on a space $\calX$ 
(typically the sample space)   
from the left 
$(g, x) \mapsto gx : \ \calG \times \calX \to \calX.$
Let $\calG x = \{gx: g \in \calG \}$ be 
the {\it orbit} containing $x \in \calX,$ and let 
$\calX/\calG = \{ \calG x : x \in \cal X \}$ be the 
{\it orbit space}, i.e., the set of all orbits. 
When $\calX$ consists of a single orbit $\calX = \calG x,$ 
the action is said to be {\it transitive}.

Indicate by 
$\calG_x = \{ g \in \calG : gx = x \}$ 
the {\it isotropy subgroup} at $x \in \calX.$
When $\calG_x = \{ e \}$ for all $x\in \calX,$ 
the action is said to be {\it free,} 
where $e$ denotes the identity element of $\calG.$
In general, the isotropy subgroups at two points on a common orbit 
are conjugate to each other:
\begin{equation}
\label{eq:G_gx}
  \calG_{gx} = g \calG_x g^{-1}, \ \ \ g \in \calG, \ x \in \calX.
\end{equation}
The set of {\it left cosets} 
$g \calG_x = \{ gg': g'\in \calG_x \}, \ g \in \calG,$ is 
called the {\it left coset space} of $\calG$ modulo $\calG_x,$ 
and is denoted by 
$\calG/\calG_x = \{ g \calG_x: g \in \calG \}.$ 
The group $\calG$ acts on $\calG/\calG_x$ by 
$(g, h\calG_x) \mapsto (gh)\calG_x, \ g,h \in \calG.$
We define the {\it canonical map} $\pi: \calG \to \calG/\calG_x$ 
by $\pi(g) = g \calG_x, \ g \in \calG.$

We move on to the definitions concerning cross sections. 
A {\it cross section} is defined to be a set 
$\calZ \subset \calX$ which 
intersects each orbit $\calG x, \ x \in \calX,$ exactly once.
Therefore, $\calZ$ is in one-to-one correspondence with 
the orbit space.
We denote this correspondence by
$\iota_\calZ : \calX/\calG \rightarrow \calZ,$ i.e., 
$\iota_\calZ(\calG x)=z,$ where $z$ is the unique point in 
$\calG x \cap \calZ.$
A cross section $\calZ$ is called a {\it global cross section} if
the isotropy subgroups are common at all points of $\calZ: \ 
\calG_z = \calG_0,$ say, for all $z \in \calZ.$ 
Of course, there always exists a cross section, but this 
is not always the case with a global cross section.

Suppose there does exist a global cross section $\calZ$ 
with the common isotropy subgroup $\calG_0.$ 
It is well-known and easy to see that 
in this case, we have 
the following one-to-one correspondence, called the {\it orbital decomposition}:
\begin{eqnarray}
\label{eq:X<->YxZ}
\calX & \leftrightarrow & \calY \times \calZ, \\
x & \leftrightarrow & (y, z), \ \ \ x = gz, \ \ y = \pi(g), 
\ \ g \in \calG, \nonumber
\end{eqnarray}
where $\calY=\calG/\calG_0$ is the left coset space modulo $\calG_0$.
It might help to regard $y$ as the {\it coordinate along the orbit} or 
the {\it within-orbit coordinate}, and 
$z$ as the {\it orbit index}.
In the orbital decomposition, we can think of $y$ and $z$ as 
functions $y=y(x)$ and $z=z(x)$ of $x.$
If $x \leftrightarrow (y,z),$ then 
$gx \leftrightarrow (gy, z), \ g\in\calG.$ 
Therefore, $y(x)$ is equivariant and $z(x)$ is invariant:
\[
 y(gx)=gy(x), \qquad z(gx)=z(x), \qquad  g \in\calG, \ x \in \calX.
\]
Thus, the coordinate along the orbit $y \in \calY$ 
is called the {\it equivariant part}, 
and the orbit index $z \in \calZ$ is called the {\it invariant part}.

From now on we assume that a global cross section $\calZ$ exists.  
We note that our results can be applied to the case of non-existence 
of a global cross section, by using the notion of orbit types.   
We discuss this point in Appendix \ref{app:orbit-types}.

We end this subsection by giving 
two simple examples 
of the orbital decompositions.

First, consider the rotation group 
\[
\calG=SO(2)=\left\{ 
\begin{pmatrix}
\cos \theta & -\sin \theta \\ 
\sin \theta & \cos \theta
\end{pmatrix}
: 0 \le \theta < 2\pi \right\}
\]  
acting on 
$\calX=\bbR^2 - \{ {\bf 0} \}
=\{ x=(x_1, x_2) \in \bbR^2 : x\ne {\bf 0} \}$ 
as 
\[ 
\left( 
\begin{pmatrix}
\cos \theta & -\sin \theta \\ 
\sin \theta & \cos \theta
\end{pmatrix},
\begin{pmatrix} x_1 \\ x_2 \end{pmatrix} 
\right)
\mapsto 
\begin{pmatrix}
\cos \theta & -\sin \theta \\ 
\sin \theta & \cos \theta
\end{pmatrix}
\begin{pmatrix} x_1 \\ x_2 \end{pmatrix} 
=
\begin{pmatrix}
(\cos \theta)x_1-(\sin \theta)x_2 \\ 
(\sin \theta)x_1+(\cos \theta)x_2
\end{pmatrix}. 
\] 
In this case, $\calG_0$ is trivial, so $\calY=\calG.$ 
The orbit containing $x=(x_1, x_2)$ is 
the circle with center ${\bf 0}$ %
and radius $\| x \| =\sqrt{x_1^2+x_2^2}.$ 
Therefore, any subset of $\bbR^2-\{ {\bf 0} \}$ intersecting each 
concentric circle with center ${\bf 0}$ exactly once can serve as a 
cross section $\calZ.$ 
We take the positive part of the $x_1$-axis as a standard cross section: 
$\calZ=\{ (x_1, 0) : x_1>0 \}.$ 
Then the orbital decomposition of $x=(x_1, x_2)$ 
can be written as 
\[ 
\begin{pmatrix} x_1 \\ x_2 \end{pmatrix}=
\begin{pmatrix}
\frac{x_1}{\| x \|} & -\frac{x_2}{\| x \|} \\ 
\frac{x_2}{\| x \|} & \frac{x_1}{\| x \|} 
\end{pmatrix}
\begin{pmatrix}
\| x \| \\ 
0
\end{pmatrix}. 
\] 
Thus we see that the equivariant part 
\[
\begin{pmatrix}
\frac{x_1}{\| x \|} & -\frac{x_2}{\| x \|} \\ 
\frac{x_2}{\| x \|} & \frac{x_1}{\| x \|} 
\end{pmatrix}
\]
can be labeled by the angle 
$\theta=\cos^{-1}(x_1/\| x \|)=\sin^{-1}(x_2/\| x \|)$ 
(the unique $\theta$ such that 
$\cos \theta =x_1/\| x \|, \ \sin \theta = x_2/\| x \|$), or 
the direction of $x.$  
The invariant part, on the other hand, can be indexed by 
the length $\| x \|.$

Now we move on to the next example. 
Let $\calG=\bbR_+^*,$ the multiplicative group of positive real numbers. 
Then $\calG$ acts freely on $\calX=\bbR^2-\{ {\bf 0} \}$ by  
$(g, (x_1, x_2))\mapsto (gx_1, gx_2).$  
The orbit containing $x=(x_1, x_2)$ 
is the ray emanating from the origin in the direction of 
$x.$ 
Hence, the cross sections are (boundaries of) ``star-shaped'' sets 
(see Section \ref{sec:star-shaped} for the precise definition). 
Here we take 
$\calZ=\bbS^1
=\{ x=(x_1, x_2) : \| x \| =1 \},$
the unit circle. 
For this cross section, $x=(x_1, x_2)$ can be 
factored as 
\[
\begin{pmatrix} x_1 \\ x_2 \end{pmatrix}=
\| x \|
\begin{pmatrix} \frac{x_1}{\| x \|} \\ 
\frac{x_2}{\| x \|} \end{pmatrix},
\] 
so the equivariant part is the length and the invariant part is the direction.

\subsection{Properties of global cross sections}

For an arbitrary (not necessarily global)
cross section $\calZ,$ we can see that 
$g \calZ = \{ gz: z \in \calZ \}$ 
is again a cross section for each $g \in \calG.$
We call $g \calZ$ a cross section {\it proportional} to $\calZ.$
Since $\calZ$ meets each orbit, we have
\begin{equation}
\label{partition-by-proportion}
\calX = \bigcup_{g \in \calG } g \calZ.
\end{equation}
We are interested in the case where (\ref{partition-by-proportion})
gives a partition of $\calX,$ that is, 
\begin{equation*}
\label{arrow}
  g_1 \calZ \cap g_2 \calZ \ne \emptyset 
  \quad \Rightarrow \quad g_1 \calZ = g_2 \calZ
\end{equation*}
for $g_1, g_2 \in \calG$.
The following proposition shows that a necessary and sufficient 
condition for %
(\ref{partition-by-proportion}) to give a partition 
of $\calX$ is that $\calZ$ be a global cross section.

\begin{prop} \quad
\label{prop:partition}
A cross section $\calZ$ is global if and only if 
$\calX = \bigcup_{g \in \calG } g \calZ$ 
gives a partition of $\calX.$
\end{prop}

\Proof
Suppose that 
$\calX = \bigcup_{g \in \calG} g \calZ$ 
gives a partition of 
$\calX.$
Let $z_1$ and $z_2$ be two arbitrary points of $\calZ.$
Let $g\in \calG_{z_1}.$ 
Then $gz_1 = z_1 \in g \calZ \cap \calZ \ne \emptyset,$ and 
hence $g \calZ = \calZ.$ 
Thus there exists a $z \in \calZ$ such that $gz_2 = z.$
But since $\calZ$ is a cross section, we have $z_2 = z$ and hence
$gz_2 = z_2.$
This observation shows that $g \in \calG_{z_1}$ implies
$g \in \calG_{z_2}.$
By interchanging the roles of $z_1$ and $z_2,$ 
we see that the converse is true as well and 
thus $\calG_{z_1} = \calG_{z_2}.$
Hence, $\calZ$ is global.

Conversely, suppose that $\calZ$ is global, and let $\calG_0$
be the common isotropy subgroup. 
Suppose $g_1 \calZ \cap g_2 \calZ \ne \emptyset$ for 
$g_1, g_2 \in \calG.$
Then, there exist %
$z_1, z_2 \in \calZ$ such that
$g_1z_1 = g_2z_2.$
Since $\calZ$ is a cross section, we have $z_1 = z_2$ and thus
$g_1z_1 = g_2z_1.$
Therefore, $g_1^{-1}g_2 \in \calG_0$ and %
$g_1z = g_2z$ for all $z \in \calZ.$
Thus we obtain $g_1 \calZ = g_2 \calZ.$
\QED

For a global cross section $\calZ,$ we call %
$\{ g \calZ : g \in \calG \}$  the {\it family of proportional
global cross sections}.

In the preceding discussions, a global cross section $\calZ$ was given 
first and the equivariant function $y$ was induced by the orbital
decomposition with respect to $\calZ.$
Conversely, we can construct a global cross section from a given
equivariant function in the following way. 
The proof of the following proposition is not difficult and is omitted.

\begin{prop} \quad
\label{prop:g.c.s. as an inverse image}
Let a group $\calG$ act on a space $\calY$ 
as well as on $\calX,$ 
and let $\tilde{y}: \calX \to \calY$ be an equivariant function. 
Suppose that the action of $\calG$ on $\calY$ is transitive 
and that %
$\tilde{y}$ satisfies the following condition:
\begin{equation}
\label{eq:y(x)=y(gx)}
  \tilde{y}(x)=\tilde{y}(gx) \ \Leftrightarrow \ x=gx
\end{equation}
for $g \in \calG$ and $x \in \calX.$
Then, the inverse image $\tilde{y}^{-1}( \{ y_0 \} ) \subset \calX$ 
of each $y_0 \in \calY$ is a global cross section. 
Moreover, global cross sections 
$\tilde{y}^{-1}(\{ y \}), \ y \in \calY,$ are all proportional 
to one another. %

\end{prop}

\begin{rem} \quad 
\label{rem:Y=G/G_0}
\begin{enumerate}
\item 
When $\calY$ {\it is} the coset space $\calY=\calG/\calG_0$ 
modulo a subgroup $\calG_0,$ 
the common isotropy subgroup of the global cross section 
$\tilde{y}^{-1}(\{ y_0 \})$ with $y_0=\calG_0$ 
coincides with $\calG_0.$ 
Furthermore, in this case $\tilde{y}(x)$ is the equivariant part of 
$x \in \calX$ with respect to this global cross section. 
Note that since the action of $\calG$ on $\calY$ is assumed to be
transitive, %
we may assume without loss of generality that $\calY$ is a
coset space. 
\item If the action of $\calG$ on $\calY$ is free, then 
$\tilde{y}(x)=\tilde{y}(gx)=g\tilde{y}(x)$ implies $g=e,$ 
so condition (\ref{eq:y(x)=y(gx)}) is 
satisfied for any equivariant function $\tilde{y}.$ 
In particular, when %
$\calY=\calG,$   
condition (\ref{eq:y(x)=y(gx)}) is automatically satisfied and 
the action of $\calG$ on $\calX$ is free, as long as an equivariant 
function $\tilde{y}: \calX \to \calG$ exists. 
\end{enumerate}
\end{rem}

When $\calY$ is a coset space $\calG/\calG_0,$ we will call 
the global cross section $\tilde{y}^{-1}(\{ \calG_0 \})$ the 
{\it unit global cross section.}

We now consider the variety of global cross sections.  {}From a
given %
global cross section $\calZ,$ we can construct a general cross section
$\calZ'$ by moving the points of $\calZ$ within their orbits.  
For example, in the case of star-shaped distributions discussed in
Section \ref{sec:star-shaped}, we consider transforming an 
ellipse $\calZ$ centered at the origin to the unit circle $\calZ'$  
by the transformation $x \mapsto x/\Vert x\Vert.$  
In the case of non-free actions, 
for $\calZ'$ to be global, i.e., for the
isotropy subgroups to be the same on the whole of $\calZ',$ 
movements of the points within the orbits have to be made 
subject to some restriction.  
Let
\[ 
  \calN=\{ g \in \calG: g \calG_0 g^{-1}=\calG_0 \}
\]
denote the normalizer of the common isotropy subgroup $\calG_0$ of a
global cross section $\calZ.$ 
Note that $\calG_0$ is a normal subgroup of $\calN$ so that we can 
think of the factor group $\calM=\calN/\calG_0$ 
(Appendix \ref{app:variety}).  
We can characterize a general global cross section in terms of 
the normalizer $\calN.$

\begin{theorem} \quad
\label{th:variety}
Let $\calZ$ be a global cross section with the common 
isotropy subgroup $\calG_0.$
Then $\calZ' \subset \calX$ is a global cross section if and only if 
it can be written as 
\begin{equation}
\label{eq:g0nz}
\calZ'=\{ g_0 n_z z: z \in \calZ \}
\end{equation}
for some $g_0 \in \calG$ and $n_z \in \calN, \ z \in \calZ.$ 
\end{theorem}

The proof is given in Appendix \ref{app:variety} 
(Corollary \ref{co:characterization}). 
As can be seen there, 
it is easy to show that $\calZ'=\{ g_0 n_z z: z \in \calZ \}$ is a
global cross section.  
The point is the proof of the converse.
Characterization of a general global cross section, including the 
proof of the converse and the question of
the uniqueness of $n_z$ in representation (\ref{eq:g0nz}),  
is fully discussed in Appendix \ref{app:variety}. %
Note that $g_0$ in (\ref{eq:g0nz}) is not essential, since
$\calZ' =\{g_0 n_z z: z \in \calZ \}$ and 
$g_0^{-1}\calZ' =\{n_z z: z \in \calZ \}$ are proportional 
and thus induce the same family of proportional global cross sections.
\begin{rem} \quad 
\label{rem:variety(free)}
When the action is free, we have $\calN=\calG$ so that  
(\ref{eq:g0nz}) becomes 
\begin{equation}
\label{eq:Z'=g_zz}
\calZ'=\{ g_z z : z \in \calZ \}
\end{equation}
for some $g_z \in \calG, \ z \in \calZ.$
\end{rem}

We finish this subsection by explicitly writing down
how the equivariant part transforms by the construction of 
a general global cross section in (\ref{eq:g0nz}).
Let $x \leftrightarrow (y,z)$ be the orbital decomposition
with respect to the %
global cross section $\calZ$ 
with the common isotropy subgroup $\calG_0,$ and 
let $ x \leftrightarrow (y', z')$ be the orbital decomposition
with respect to 
the $\calZ'$ in (\ref{eq:g0nz}). 
This $\calZ'$ has the isotropy subgroup 
$%
\calG_0'=g_0 \calG_0 g_0^{-1}.$
Now the equivariant part based on $\calZ'$ is given 
as follows: 

\begin{prop} \quad
\label{prop:transformation of equivariant part}
Let $\calZ$ be a global cross section with the common 
isotropy subgroup $\calG_0.$
Moreover, let $\calZ'$ be as in (\ref{eq:g0nz}), and 
$x \leftrightarrow (y',z')$ the orbital decomposition 
with respect to $\calZ'.$  
Write $x \in \calX$ as $x = gz = g'z'$ with 
$z \in \calZ, \ z' \in \calZ'$ and $g, g'\in \calG.$
Then we have 
$
y' = y n_z^{-1} g_0^{-1},
$
where $y = g \calG_0$ and 
$y' = g' \calG'_{0}.$
\end{prop}
\Proof 
We can write $x=gz$ as $x = g n_z^{-1} g_0^{-1} z'$ in terms of 
$z'=g_0 n_z z.$
This implies $y' = (g n_z^{-1} g_0^{-1})(g_0 \calG_0 g_0^{-1})
=g \calG_0 n_z^{-1} g_0^{-1}
=y n_z^{-1} g_0^{-1}.$
\QED
\begin{rem} \quad 
\label{rem:transformation of equivariant part (free)}
When the action is free 
and a general cross section is given by (\ref{eq:Z'=g_zz}), 
the equivariant part transforms as  
$g'=g g_z^{-1},$ where $g$ and $g'$ are the equivariant parts 
with respect to $\calZ$ and $\calZ',$ respectively. 
\end{rem}

\section{Decomposable distributions}
\label{sec:decomposable-distribution}

In this section we define a class of distributions called decomposable
distributions and study some distributions induced by them.  
The general discussion here is applied to particular cases in the next two sections. 
Especially, 
an extension of elliptically contoured distributions called star-shaped
distributions is discussed in Section \ref{sec:star-shaped}. 

\subsection{Assumptions and the definition}

In order to make distributional arguments, we need to make 
topological and measure-theoretic 
assumptions.
In this paper, measurability of topological spaces refers to %
Borel measurability.

\begin{assumption} \quad
\label{assumption:regularity1}
\begin{enumerate}
\setlength{\itemsep}{0pt} 
 \item $\calX$ is a locally compact Hausdorff space.
 \item $\calG$ is a second countable, locally compact Hausdorff
 topological group acting continuously on $\calX.$
 \item $\calG_0$ is compact.
\item Global cross section $\calZ$ is locally compact, and the bijection
$x \leftrightarrow (y, z)$ with respect to $\calZ$ 
is bimeasurable, where
the topology on $\calZ$ is the relative topology of $\calZ$ as
a subset of $\calX.$
\end{enumerate}
\end{assumption}

We agree that a quotient space receives the quotient topology 
when regarded as a topological space.  
This applies to the coset space $\calY = \calG/\calG_0$ as well as to 
the orbit space $\calX/\calG.$ 
Because of 2 of Assumption \ref{assumption:regularity1}, 
there exists a left Haar measure $\mu_{\calG}$ on $\calG,$ 
which is unique up to a multiplicative constant.  

We consider densities 
with respect to a dominating measure $\lambda$ on $\calX$ which is 
relatively invariant with multiplier $\chi:$ 
\begin{equation*}
\label{eq:relative-invariant-measure}
\lambda(d(gx))=\chi(g)\lambda(dx), %
\ g \in \calG.
\end{equation*}
Note that the relative invariance of $\lambda$ only determines its 
behavior within each orbit so  
that for any nonnegative $f_{\calZ}(z(x)),$ 
\begin{equation}
\label{eq:relative-invariance}
\tilde{\lambda}(dx) = f_{ \calZ }(z(x))\lambda(dx)
\end{equation}
is again a relatively invariant measure with the same multiplier 
$\chi$ as $\lambda.$

We are now in a position to define the decomposable distributions.
\begin{definition} \quad
A distribution on $\calX$ is said to be decomposable with 
respect to a global cross section $\calZ$ if it is of the form
\[
  f(x)\lambda(dx) = f_{ \calY }(y(x))f_{ \calZ }(z(x))\lambda(dx).
\]
In particular, it is said to be cross-sectionally contoured if 
$f_{ \calZ }(z)$ is constant. 
In contrast, it is said to be orbitally contoured if 
$f_{ \calY }(y)$ is constant. 
\end{definition}

Obviously, a distribution $f(x)\lambda(dx)$ is
cross-sectionally contoured with respect to $\calZ$ 
if and only if $f(x)$ is constant on each
proportional global cross section $g\calZ, \ g \in \calG.$
Similarly, $f(x)\lambda(dx)$ is orbitally contoured if and only if $f(x)$ is
constant on each orbit $\calG x, \ x \in \calX.$
Before examining the distributions of the invariant and equivariant
parts, %
we observe the following two points.

First, a decomposable distribution
$
  f_{ \calY }(y(x)) f_{ \calZ }(z(x)) \lambda(dx)
  = f_{ \calY }(y(x))\tilde{\lambda}(dx) 
$
can always 
be thought of as a cross-sectionally contoured distribution in view of 
(\ref{eq:relative-invariance}).

Next, we can take various global cross sections, in addition to
``standard'' ones like the unit sphere. %
This enables us to consider the cross-sectionally
contoured distributions associated with a variety of global
cross sections.
In contrast, once an action is given, there is no room for
choosing the orbits; 
the orbits are determined by the action in question, and usually
those orbits are familiar subsets of $\calX.$
Hence, we can not produce the orbitally contoured distributions
based on the orbits which are unfamiliar subsets of $\calX.$ 

For these reasons, we will be concerned %
with the cross-sectionally contoured distributions from now on.

\subsection{Distributions of invariant and equivariant parts}

First, we confirm the independence of invariant and 
equivariant parts. 
This corresponds to the independence of ``direction'' and 
``length'' in elliptically contoured distributions.
Thanks to the assumption that $\calG_0$ is compact, 
we have the induced measure 
$\mu_{ \calY } = \pi(\mu_{ \calG }) = \mu_{ \calG } \pi^{-1}$ 
on $\calY$ (Proposition 2.3.5 and Corollary 7.4.4 of \cite{Wijsman90}).
Also, by the same assumption we can define 
$\bar{\chi}(y), \ y \in \calY,$ by 
$\bar{\chi}(y) = \chi(g)$ with $g \in \pi^{-1}(\{ y \}),$
where $\chi$ is the multiplier of $\lambda.$ 
With some abuse of notation,
we will write $\chi(y)$ for $\bar{\chi}(y).$

Now, %
$\lambda(dx)$ is factored as   
\begin{equation}
\label{eq:decomposition of lambda}
\lambda(dx) = \chi(y)\mu_{ \calY }(dy) \nu_{ \calZ }(dz)
\end{equation}
(Theorem 7.5.1 of \cite{Wijsman90}, 
Theorem 10.1.2 of \cite{Farrell}). 
By changing the weights of the orbit as (\ref{eq:relative-invariance}) if
necessary, from now on we assume that $\nu_{ \calZ }(dz)$ is 
(standardized to be) a probability measure on $\calZ$. 
The following theorem is an
immediate consequence of the factorization of $\lambda(dx)$ in
(\ref{eq:decomposition of lambda}).

\begin{theorem} \quad
\label{th:distributions of y and z}
Suppose that $x$ is distributed according to a cross-sectionally
contoured distribution $f_{\calY}(y(x))\lambda(dx).$
Then we have:
\begin{enumerate}
\item $y = y(x)$ and $z = z(x)$ are independently distributed.
\item The distribution of $z$ does not depend on $f_{ \calY }.$
\item The distribution of $y$ is 
$f_{\calY}(y)\chi(y)\mu_{ \calY }(dy).$
\end{enumerate}
\end{theorem}

\subsection{Distributions generated via two global cross sections}

Next, we investigate distributions generated by considering 
two global cross sections at a time.
The relation between two global cross sections was given in 
Theorem \ref{th:variety}.
Let $\calZ$ and $\calZ'$ be two global cross sections.
By choosing an appropriate global cross section from the family of
proportional global cross sections, we assume without essential loss 
of generality that %
the common isotropy subgroups for $\calZ$ and %
$\calZ'$ are the same:
$\calG_{z} = \calG_{z'} = \calG_0, \ z \in \calZ, \ z' \in \calZ'.$

The invariant and equivariant parts $z=z(x), \ y=y(x)$ with 
respect to $\calZ$ are given via the orbital decomposition 
(\ref{eq:X<->YxZ}) as before. 
In a similar manner, define the invariant and equivariant parts 
$z' = z'(x), \ y' = y'(x)$ with respect to $\calZ'$.
Note that by our assumption the coset spaces are the 
same, $\calY=\calG/\calG_0,$ in both cases.

Denote by $g: \calY \to \calG$ an arbitrary selection 
$g(y) \in y \subset \calG.$ 
{}From now on, we will write $g(x)$ for $g(y(x)): \ 
x = g(x)z(x), \ x \in \calX.$
Define $g'(x)$ in the same way:
$x = g'(x)z'(x), \ x \in \calX.$
Here we define the map $w: \calX \rightarrow \calX$ by 
\begin{equation*}
\label{within-orbit map}
w = w(x) = g(x) z'(x)
= g(x)g'(x)^{-1}x 
= g(x)g'(z(x))^{-1}z(x)
, \ \ \ x \in \calX.
\end{equation*}
Note that since $\calG_{z'} = \calG_0, \ z' \in \calZ'$, 
$w$ does not depend on the choice of the 
selection $g(y).$
We call $w$ the {\it within-orbit bijection},   
because $x$ and $w(x)$ are on the same orbit and we are transforming
$x$ to $w(x)$ in each orbit separately. 
The within-orbit bijection
is a basic tool for deriving a new cross-sectionally
contoured distribution from a given cross-sectionally
contoured distribution. 

\begin{theorem} \quad
\label{th:distribution of w}
Suppose that $x$ is distributed according to a cross-sectionally
contoured distribution $f_{\calY}(y(x))\lambda(dx).$ 
Then the distribution of $w = w(x)$ is
\[
f_{\calY}(y'(w))\chi(g(w)^{-1}g'(w))
\Delta^{ \calG }(g(w)^{-1}g'(w))\lambda(dw),
\]
where $\Delta^{ \calG }$ is the right-hand modulus of $\calG:$
$
\mu_{\calG}(d(gg_1)) = 
\Delta^{ \calG }(g_1)\mu_{\calG}(dg), \ %
g_1 \in \calG.
$
\end{theorem}

\Proof
We regard $w=w(x)$ as a function of $y=y(x)$ and $z=z(x):$
$w = w(x) = w(y,z).$
Noting that the integration over $\calY$ can be carried out by the
integration over $\calG,$ we have
for an arbitrary measurable subset $B \subset \calX$ that
\begin{eqnarray}
\label{eq:P(wB)}
P(w \in B)
&=& \int_{ \calX } I_B(w(x))f_{\calY}(y(x))\lambda(dx) \notag \\
&=& \int_{ \calZ } \int_{ \calY} I_B(w(y,z))f_{\calY}(y)\chi(y)
\mu_{ \calY }(dy)\nu_{ \calZ }(dz) \notag \\
&=& \int_{ \calZ} \int_{ \calG } I_B(gg'(z)^{-1}z)
\hat{f}_{ \calY }(g) \chi(g) \mu_{ \calG }(dg) \nu_{ \calZ }(dz),
\end{eqnarray}
where $\hat{f}_{\calY} = f_{\calY} \circ \pi$ and $I_B$ is the indicator 
function of $B.$
Let $\calG_B(z) = \{ g \in \calG: gz \in B \}, \ z \in \calZ,$ 
and $\calZ_B = \{ z(x) \in \calZ: x \in B \}.$ 
Then $I_B(x)=I_{\calG_B(z)}(g) \cdot I_{\calZ_B}(z)$ for $x=gz,$ and 
we can write (\ref{eq:P(wB)}) as  
\begin{eqnarray*}
&& \int_{\calZ} I_{\calZ_B}(z) \int_{\calG} I_{\calG_B(z)}(gg'(z)^{-1}) 
\hat{f}_{ \calY }(g) \chi(g) \mu_{ \calG }(dg) \nu_{ \calZ }(dz) \\ 
&& \ \ \ = \int_{\calZ} I_{\calZ_B}(z) \chi(g'(z)) \Delta^{\calG}(g'(z))
\int_{\calG} I_{\calG_B(z)}(g) \hat{f}_{\calY}(gg'(z)) 
\chi(g) \mu_{ \calG }(dg) \nu_{ \calZ }(dz) \\ 
&& \ \ \ = \int_{\calX} I_B(x) \hat{f}_{\calY}(g'(x))
\chi(g(x)^{-1}g'(x))\Delta^{\calG}(g(x)^{-1}g'(x)) \lambda(dx) \\ 
&& \ \ \ = \int_{B} f_{\calY}(y'(w))
\chi(g(w)^{-1}g'(w))\Delta^{\calG}(g(w)^{-1}g'(w)) \lambda(dw).
\end{eqnarray*}
\QED

For notational simplicity, we will write
\begin{equation*}
\label{eq:delta-g}
\Delta(g) = \chi(g)\Delta^{ \calG }(g), \ \ \ g \in \calG, 
\end{equation*}
which is a continuous homomorphism from $\calG$ to $\bbR_+^*.$ 
Because of 3 of Assumption \ref{assumption:regularity1}, 
we have $\Delta(g) = 1$ for all $g \in \calG_0,$ 
so $\Delta(g(w)^{-1}g'(w))$ does not depend on
the choice of the selections $g(w)$ and $g'(w).$

Let ${\cal E}$ be the set of all measurable  
equivariant functions  
$\tilde{y}: \calX \to \calG/\calG_0$ satisfying 
(\ref{eq:y(x)=y(gx)})
of Proposition \ref{prop:g.c.s. as an inverse image} 
and let $\tilde{g}: \calX \to \calG$ be an arbitrary selection of 
$\tilde{y}:$ $\tilde{g}(x) \in \tilde{y}(x), \ x \in \calX.$ 
By Theorem \ref{th:distribution of w},
when $f_{\calY}(y(x))$ is a density function with respect to $\lambda(dx),$
we can define a non-parametric family of distributions dominated by $\lambda:$
\begin{equation}
\label{distribution-family-2}
\left\{ f(x; h, \tilde{y}) = \frac{1}{\chi(h)} f_{\calY}(h^{-1}\tilde{y}(x)) \Delta(g(x)^{-1} \tilde{g}(x)) : 
h \in \calG, \ \tilde{y} \in {\cal E} \right\}. 
\end{equation}
Note that $\tilde{y}(x)$ is the equivariant part of $x$ with respect to 
the unit global cross section $\tilde{y}^{-1}(\{ \calG_0 \})$ 
(Remark \ref{rem:Y=G/G_0}). 
We can see that distributions in %
(\ref{distribution-family-2}) have 
cross-sectionally contoured densities 
$\chi(g)^{-1}f_{\calY}(g^{-1}\tilde{y}(x))$ 
with respect to global cross section $\tilde{y}^{-1}(\{ \calG_0 \})$ 
and %
dominating measure 
$
  \tilde{\lambda}(dx) = \Delta(g(x)^{-1}\tilde{g}(x))\lambda(dx).
$
This 
$\tilde{\lambda}$ is relatively invariant with the same 
multiplier $\chi$ as $\lambda,$ and  
$\tilde{\lambda}$ and $\lambda$ are absolutely continuous with respect to 
each other because $ 0 < \Delta(g(x)^{-1}\tilde{g}(x)) < \infty.$

Now we turn to the distribution of $z'=z'(x).$ 
Note that we may instead obtain the distribution of $z'(w)$ because 
$z'(x)=z'(w(x)),$ for $x$ and $w(x)$ are on the same orbit. 
Corresponding to the orbital decomposition with respect to
$\calZ'$,  $\lambda(dx)$ is factored as
\begin{equation*}
\label{eq:decomposition of lambda wrt Z'}
  \lambda(dx)
  = %
     \chi(y')\mu_{ \calY }(dy')
  \nu_{ \calZ' }(dz').
\end{equation*}
Here we use the same $\mu_{ \calY }$ as in 
(\ref{eq:decomposition of lambda}).  
Recall that in (\ref{eq:decomposition of lambda}) we have chosen 
$\nu_{ \calZ }(dz)$ to be a probability measure on $\calZ.$ 
Therefore, $\nu_{ \calZ' }$ is not necessarily a probability 
measure on $\calZ'.$
In terms of $\nu_{ \calZ' },$ the distribution of $z'$
is written as follows.

\begin{theorem} \quad
\label{th:distribution of z^prime}
Suppose that $x$ is distributed according to a cross-sectionally
contoured distribution $f_{\calY}(y(x))\lambda(dx).$ 
Then the distribution of $z' = z'(x)$ is
\begin{equation}
\label{eq:distribution of z' (general)}
\frac{1}{\Delta(g(z'))}\nu_{ \calZ' }(dz').
\end{equation}
In addition, $z'=z'(x)$ is independently distributed of $y=y(x).$
\end{theorem}

\Proof
We have $\Delta(g(w)^{-1}g'(w)) = \Delta(g(z'(w)))^{-1}.$ 
Writing $y' = y'(w)$ and $z' = z'(w),$
we have by Theorem \ref{th:distribution of w} that the 
distribution of $w$ is
\begin{eqnarray}
\label{eq:y'(w)z'(w)}
f_{\calY}(y'(w))\Delta(g(z'(w)))^{-1}\lambda(dw)
       &=& f_{\calY}(y')\Delta(g(z'))^{-1}
              \chi(y')\mu_{ \calY }(dy')
           \nu_{ \calZ'}(dz')  \notag \\ 
       &=& %
           f_{\calY}(y')\chi(y')\mu_{ \calY }(dy')
           \Delta(g(z'))^{-1}\nu_{\calZ'}(dz').
\end{eqnarray}
Accordingly, the distribution of $z' = z'(w)$ is
$\Delta(g(z'))^{-1}\nu_{ \calZ' }(dz').$

Since $x$ and $w=g(x)g'(x)^{-1}x$ are on the same orbit,
we have $z'(x) = z'(w)$ so that the distribution of
$z'(x)$ is the same as that of $z'(w).$
Moreover, we can see from (\ref{eq:y'(w)z'(w)}) that 
$y'(w)=g(x)g'(x)^{-1}y'(x)=y(x)$ and $z'(w)=z'(x)$ are independent. 
\QED

{}From (\ref{eq:distribution of z' (general)})  
we can construct various distributions on $\calZ'$ by appropriately
choosing %
the global cross sections $\calZ.$ 
Here we can ask the following question: 
Given a density $f(z')$ on $\calZ',$ can we find a global cross section 
$\calZ$ such that the distribution of 
$z'(x)$ when $x$ is 
distributed as a cross-sectionally contoured distribution 
$f_{\calY}(y(x))\lambda(dx)$ with respect to $\calZ$ 
coincides with $f(z')$?  
Recall that the distribution of $z'(x)$ depends only on $\calZ$ and 
not on $f_{\calY}.$
The following corollary gives the answer.

\begin{co} \quad
\label{co:given-density-on-Z'}
Let $f(z')\nu_{\calZ'}(dz')$ be a distribution on $\calZ'$ such that
$f(z')$ is almost everywhere positive on $\calZ'$ with respect
to $\nu_{\calZ'}.$  
Suppose there exists a coset 
$%
g\calN', \ g \in \calG,$ 
with respect to the normalizer $%
\calN'$ of 
$%
\calG_0'=\calG_{z'}, \ z' \in \calZ',$ such that 
\begin{equation}
\label{eq:ne-empty}
\Delta^{-1}(\{ f(z') \}) \cap g %
\calN' \ne \emptyset
\end{equation} 
for each $z' \in \calZ'$ with positive $f(z').$ 
Then there exists a global cross section $\calZ$ such that 
the distribution of $z'=z'(x)$ coincides with $f(z')\nu_{\calZ'}(dz')$ for
$x$ having an arbitrary cross-sectionally contoured distribution with
respect to $\calZ$.
\end{co}

\Proof  
For any $z' \in \calZ'$ with $f(z')>0,$  
we can choose $g(z') \in \calG$ such that 
$g(z')^{-1} \in \Delta^{-1}(\{ f(z') \}) \cap g \calN'.$
Take a cross section $\calZ=\{ g(z')^{-1}z': z' \in \calZ' \}.$ 
Then $\calZ$ is global with the common isotropy subgroup 
$g\calG_0' g^{-1}.$
Writing $z' \in \calZ'$ as $z'=g(z')\cdot g(z')^{-1}z',$ we see that 
$g(z')$ can serve as a selection of the equivariant part of $z'$ with respect to 
$\calZ.$ 
Theorem \ref{th:distribution of z^prime} implies that the density of 
$z'=z'(x)$ is $\Delta(g(z')^{-1})=f(z').$
\QED

\begin{rem} \quad
\label{rem:given-density-on-Z'}
When the action is free, we have $\calN'=\calG$ and 
condition (\ref{eq:ne-empty}) is satisfied as long as   
$\Delta$ is not identically equal to $1$
since $\Delta$ is a continuous homomorphism from $\calG$ to
$\bbR_+^*.$ 
If in addition $\calG$ is unimodular 
$\Delta^{\calG} \equiv 1$ (e.g., abelian or compact), then  
condition $\Delta = \chi \not\equiv 1$ is 
equivalent to $\lambda$ not being an invariant measure. 
\end{rem}

\section{Star-shaped distributions}
\label{sec:star-shaped}

In this section, we define star-shaped distributions in $\bbR^p$ 
and investigate their properties.
Most results presented here are easy consequences of the general 
arguments in the preceding section, 
but also included here are results which can be obtained only after %
regarding the orbits and cross sections as submanifolds of $\bbR^p.$ 

Let $\calG=\bbR_+^*$ %
and define its action on $\calX=\bbR^p - \{ {\bf 0} \}$ by
\begin{equation}
\label{eq:action(star-shaped)}
(g, (x_1, \ldots, x_p)) \mapsto (gx_1,\ldots,gx_p). 
\end{equation}
Under this action, the Lebesgue measure $dx$ is relatively invariant 
with multiplier $\chi(g) = g^p.$ 
We take $\lambda(dx)=dx$ as the dominating measure.
Note that the origin has %
Lebesgue measure zero and that omitting it 
in the sample space $\calX=\bbR^p-\{ {\bf 0} \}$ does not affect 
the discussion about the distributions in $\calX.$
By so doing, we have made the sample space have just one orbit 
type (Appendix \ref{app:orbit-types}) and %
made our action (\ref{eq:action(star-shaped)}) free.   

Since the action is free, we know from 
Remark \ref{rem:Y=G/G_0} that choosing a unit cross section $\calZ$ is   
equivalent to choosing an equivariant function from $\calX$ to $\bbR_+^*.$
Now, let $g: \calX \to \bbR_+^*$ be an equivariant function. 
We call distributions with the densities of the form 
\begin{equation}
\label{eq:star-shaped-density}
f(x)=f_{\calG}(g(x))
\end{equation}
{\it star-shaped} distributions. 
Obviously, this reduces to the elliptically contoured 
distributions when $g(x)=(x\trans\Sigma^{-1}x)^{1/2}$ with 
$\Sigma \in PD(p)$ (the set of $p \times p$ positive 
definite matrices).

The orbits under (\ref{eq:action(star-shaped)}) are rays 
emanating from the origin, so the unit cross section 
$\calZ = \{ x : g(x)=1 \}$ associated with $g$ 
is a set which meets each ray exactly once.  
Hence, 
\begin{equation}
\label{eq:UcZ}
\bigcup_{0 \le c \le 1} c\calZ
\end{equation}
contains every line segment connecting the origin and
a point on $\calZ.$ 
Namely, (\ref{eq:UcZ}) is a star-shaped set with respect to the origin.  
This is why we call the distributions with densities of the 
form (\ref{eq:star-shaped-density}) star-shaped. 
(For the term ``star-shaped,'' see also Definition 3.1 of 
\cite{Naiman-Wynn}.)

Throughout this section, we assume $x$ is distributed 
according to $f_{\calG}(g(x))dx.$

A version of the Haar measures on $\calG=\bbR_+^*$ is given 
by $g^{-1}dg.$ 
By Theorem \ref{th:distributions of y and z},  
$g=g(x)$ and $z=x/g(x)$ are independent and the 
joint distribution of $g$ and $z$ can be written as
\begin{equation*}
  \label{star-shaped-ind-1}
  \frac{1}{c_0} f_{\calG}(g) g^{p-1}dg \times \nu_{\calZ}(dz), 
\end{equation*}
where $c_0=\int_0^\infty f_{\calG}(g) g^{p-1}dg,$ and 
$\nu_{\calZ}$ is a probability measure on $\calZ.$
Note that we have taken $\mu_{\calG}(dg)=c_0^{-1}g^{-1}dg.$

For %
action (\ref{eq:action(star-shaped)}), the most standard 
cross section is the unit sphere
$\calZ'=\bbS^{p-1}=\{ x \in \calX : g'(x) = 1 \},$
where $g'(x)=\Vert x \Vert = (x\trans x)^{1/2}$ is the
usual Euclidean length of $x \in \bbR^p-\{ {\bf 0}\}.$ 
Now, $dx$ obviously  factors as 
\begin{equation*}
\label{eq:dx}
   dx = {g'}^{p-1} dg' \, dz' 
      = \frac{1}{c_0} {g'}^{p-1}dg' \times c_0 dz',
\end{equation*}
where $dz'$ is the volume element of $\bbS^{p-1}.$ 
Since $\calG=\bbR_+^*$ is abelian, we have 
$\Delta(g)=\chi(g)=g^p.$
Thus, by Theorem \ref{th:distribution of z^prime} 
the distribution of the direction vector 
$z'=x/\Vert x \Vert$ is obtained as
\begin{equation}
\label{eq:distribution of z'}
   c_0\, g(z')^{-p} dz', 
\end{equation}
from which another expression of $c_0$ can be given: 
$c_0 = 1 \ /\int_{\bbS^{p-1}}g(z')^{-p} dz'.$

When (\ref{eq:star-shaped-density}) is an elliptically 
contoured density, %
(\ref{eq:distribution of z'}) 
becomes 
$c_0\, ({z'}\trans\Sigma^{-1}z')^{-p/2} dz'.$
Normalizing constant $c_0,$ being 
independent of the choice of $f_{\calG}( \cdot ),$ 
can be obtained by considering
the particular case of normality
$f_{\calG}(g)=(2\pi)^{-p/2}(\det\Sigma)^{-1/2}\exp(-g^2/2)$ 
as
\[
c_0 = 
(2\pi)^{-\frac{p}{2}}(\det\Sigma)^{-\frac{1}{2}}
\int_0^\infty e^{-\frac{g^2}{2}} g^{p-1}dg 
=\omega_p^{-1}(\det\Sigma)^{-\frac{1}{2}}, 
\]
where 
$\omega_p=2\pi^{p/2}/\Gamma(p/2)=\int_{\bbS^{p-1}} dz'$ is the total 
volume of $\bbS^{p-1}.$ 
This distribution %
\begin{equation}
\label{eq:angular Gaussian}
\frac{1}{\omega_p(\det\Sigma)^{\frac{1}{2}}}
({z'}\trans\Sigma^{-1}z')^{-\frac{p}{2}} dz'
\end{equation}
is derived in Section 3.6 of \cite{Watson}.  
Our (\ref{eq:distribution of z'}) is a generalization of 
(\ref{eq:angular Gaussian}) to the case of an arbitrary (i.e., not necessarily 
elliptically contoured) star-shaped density.

\begin{rem} \quad
Distribution (\ref{eq:angular Gaussian}) has been studied  
in several parts of the literature.
\cite{Watson} also notes that (\ref{eq:angular Gaussian}) can be 
thought of as a special case of the angular Gaussian distribution, 
and discusses some of its properties.
Several arguments about statistical inferences based on this model
are given in \cite{Tyler87}.
See also Sections 9.4.4, 10.3.5 and 10.7.1 of \cite{Mardia-Jupp}.
The special case $p=2$ is treated in \cite{Kent-Tyler88}
and Section 3.5.6 of \cite{Mardia-Jupp}.

The distribution (\ref{eq:angular Gaussian}) of
$z'=x/\Vert x \Vert,$ as well as that of
${z'}\trans \Sigma^{-1}z' = x\trans \Sigma^{-1}x/x\trans x,$ 
plays an important role in null robust testing problems.
See, e.g., \cite{Kariya-Eaton} and \cite{King}.
\end{rem}

We now investigate star-shaped distributions more closely 
by viewing the orbits and cross sections as submanifolds of 
$\bbR^p.$ 
We make the additional assumption that
$g(x)$ is piecewise of class $C^1.$  

Fix $z_0 \in \calZ$ and  call $M_C(z_0)= \calZ$ 
the {\it cross section manifold} 
and $M_O(z_0)=\{ u_1 z_0 : u_1 > 0\}$ 
the {\it orbit manifold} through $z_0$.
The tangent vector of $M_O(z_0)$ at $z_0$ is 
$v_1 = z_0$. 
Choose local coordinates $u_2, \ldots, u_p$ of 
$M_C(z_0)=\calZ$ 
such that 
$v_j = \frac{\partial}{\partial u_j}
z(0, \ldots, 0, u_j, 0, \ldots, 0)\big|_{u_j=0}$, $j=2, \ldots, p,$
are orthonormal vectors.  
Then $du_2 \cdots du_p$ is the volume element of $M_C(z_0)$ at $z_0.$ 
Writing $x=u_1z(u_2,\ldots,u_p),$ we see that
\begin{equation}
   \label{star-shaped-jacobian-1}
  dx = | \det(v_1, \ldots, v_p) | \times du_1 \times du_2 \cdots du_p,
\end{equation}
where $(v_1, \ldots, v_p)$ denotes the matrix consisting 
of columns $v_1, \ldots, v_p.$ 

Let $n_{z_0}$ be the unit normal vector of $\calZ$ at $z_0$ 
pointing outward of the star-shaped set $\bigcup_{0 \le c \le 1} cZ.$   
Write $v_1=z_0$ as a linear combination of the orthonormal vectors 
$n_{z_0}, v_2, \ldots, v_p$ as 
$z_0 = a_1 n_{z_0} + a_2 v_2+ \cdots + a_p v_p.$ 
Then 
$| \det(v_1,\ldots,v_p) | = a_1 = \langle z_0,n_{z_0}\rangle=z_0\trans n_{z_0},$ 
and (\ref{star-shaped-jacobian-1}) is written as 
$dx =  du_1 \times du_2 \cdots du_p 
\times \langle z_0, n_{z_0} \rangle.$  For the rest of this section
$\langle \cdot , \cdot\rangle$ denotes the standard inner product of
${\mathbb R}^p$.
Rewrite this further as
\begin{equation}
   \label{star-shaped-jacobian-3}
  dx = \Vert z_0 \Vert du_1 \times du_2 \cdots du_p \times
   \left\langle \frac{z_0}{\Vert z_0 \Vert}, n_{z_0} \right\rangle.
\end{equation}
Note that the first term $\Vert z_0 \Vert du_1 = \sqrt{\langle v_1, v_1 \rangle} du_1$ 
in (\ref{star-shaped-jacobian-3}) 
is the volume element of $M_O(z_0)$ around $z_0.$
The second term $du_2 \cdots du_p$ is the volume element of $M_C(z_0)$ as mentioned 
above. 
Concerning the third term, let $\theta$ denote the angle between 
$z_0$ and $T_{z_0}(M_C),$ where $T_{z_0}(M_C)$ 
stands for the tangent space of $M_C(z_0)$  at $z_0$.
Then $\pi/2 -  \theta$ is the angle between $z_0$ and $n_{z_0},$ 
and the third term in (\ref{star-shaped-jacobian-3}) can be 
written as $\left\langle z_0/\Vert z_0 \Vert, n_{z_0} \right\rangle = \sin \theta.$
Therefore, (\ref{star-shaped-jacobian-3}) means that $dx$ can be 
factored into the volume elements of $M_O(z_0)$ and $M_C(z_0)$ and the 
sine of the angle between $T_{z_0}(M_O)$ and $T_{z_0}(M_C).$

We note in passing that the unit normal vector $n_{z_0}$ coincides 
with the normalized gradient of $g(x),$ i.e., 
$n_{z_0} = \nabla g(z_0)/\Vert \nabla g(z_0) \Vert.$ 
We also note the following fact.  
Let $H_{z_0} = z_0 + T_{z_0}(M_C)$ be the tangent hyperplane 
of $\calZ$ at $z_0.$  
Then
\begin{equation*}
    \label{distance-from-origin-to-tangent-hyperplane}
    \langle z_0,n_{z_0} \rangle = 
      \mbox{Euclidean distance from the origin to $H_{z_0},$}
\end{equation*}
which is the {\it support function} at $z_0$ 
(Section 8.1 of \cite{Flanders}).

Now consider the translation by $g \in \calG=\bbR_+^*$ 
from $x=z_0$ to $x=gz_0.$ 
Since this translation is just the scale change, 
its effect is straightforward. 
The volume element of the orbit manifold $M_O(z_0)$ is 
multiplied by $g,$ 
and the volume element of the cross section manifold $M_C(z_0)$ is 
multiplied by $g^{p-1},$ with $p-1$ being the dimensionality 
of $M_C(z_0).$  
Furthermore, the angle between these two manifolds remains unchanged  
under the translation. 
Therefore, around $x=gz_0$ the volume element $dx$ is 
\begin{eqnarray*}
dx 
&=& \Vert z_0 \Vert dg \times g^{p-1} dz \times 
\left\langle \frac{z_0}{\Vert z_0 \Vert}, n_{z_0} \right\rangle \\ 
&=& \frac{1}{c_0} g^{p-1}dg \times 
c_0 \, \langle z_0, n_{z_0} \rangle \,dz, 
\end{eqnarray*}
where $dg=g \,du_1$ is the volume element of 
$\calG = \bbR_+^*$ around $g \in \calG,$ 
and $dz$ is the volume element of $\calZ=M_C(z_0).$
Therefore, the distribution $\nu_{\calZ}$ of $z$ can be expressed as  
$\nu_{\calZ}(dz)=c_0 \, \langle z_0, n_{z_0} \rangle \,dz.$ 

\bigskip

We list some examples of star-shaped distributions.
\begin{itemize}
\setlength{\itemsep}{0pt}
\item[(a)] Elliptically contoured distribution: 
When $g(x)=(x\trans\Sigma^{-1}x)^{1/2},$ 
we have $\langle z, n_{z} \rangle 
= \langle z, \Sigma^{-1}z/\Vert \Sigma^{-1}z \Vert \rangle 
= (z\trans \Sigma^{-2} z)^{-1/2}$ for $z \in \calZ=\{ x: g(x)=1 \}.$ 
So in this case, $\nu_{\calZ}(dz)$ has density 
$\omega_p^{-1}(\det\Sigma)^{-1/2}(z\trans \Sigma^{-2} z)^{-1/2}$ with 
respect to the volume element of the ellipsoid 
$\{ z\in \bbR^p: z\trans \Sigma^{-1}z =1 \}.$ 
\item[(b)] ``Hypercube distribution'': Take $g(x)=
\max(|x_1|,\ldots,|x_p|), \ x=(x_1,\ldots,x_p).$ 
Then the unit cross section $\calZ$ is the surface of 
the hypercube $C_p$ in $\bbR^p$ 
and we have $\langle z, n_z \rangle=1$ on the relative interiors of 
the facets of $C_p.$  
Hence $\nu_{\calZ}(dz)$ is the uniform distribution on 
$\partial C_p.$ 
Constant value of the density is $c_0=1/{\rm Vol}_{p-1}(\partial C_p)
=1/(2p \times 2^{p-1})=1/(2^p p).$ 
\item[(c)]
 ``Crosspolytope distribution,''  
 also known as $\ell_1$\kmms-norm symmetric distribution \\ 
 (\cite{Fang87}, \cite{Fang-Anderson}, \cite{Fang-Kotz-Ng}): 
Let $g(x)=|x_1|+\ldots+|x_p|, 
\ x=(x_1,\ldots,x_p).$
The associated unit cross section $\calZ$ is the surface of 
the crosspolytope $C_p^{\Delta},$ which is polar to $C_p$ 
(Chapter 0 of \cite{Ziegler}). 
Since $\langle z, n_z \rangle,$ the distances of the facets 
of $C_p^{\Delta}$ from the origin, are constant $(=1/\sqrt{p})$ by 
symmetry, $\nu_{\calZ}(dz)$ has constant density 
$c_0/\sqrt{p}=1/{\rm Vol}_{p-1}(\partial C_p^{\Delta})
=1/\{ 2^p \times \sqrt{p}/(p-1)!\}=(p-1)!/(2^p \sqrt{p}).$ 
\item[(d)] Take $\calZ$ to be the surface of 
a $p$\kmms-dimensional polytope $P \ ({\bf 0}\in {\rm int}(P))$ 
whose facets are not equidistant from the origin.
Then we obtain a non-uniform distribution of $z$ on $\partial P.$ 
\end{itemize}

We now summarize our results in this section in the following theorem.

\begin{theorem} \quad
\label{th:star-shaped}
Suppose the distribution of $x \in \bbR^p-\{ {\bf 0} \}$ has  
a star-shaped density $f_{\calG}(g(x))$ with respect to $dx.$  
Then $g=g(x)$ and $z=x/g(x)$ are independent and
the joint distribution of $g$ and $z$ is written as
\begin{equation}
\label{eq:joint-star}
\frac{1}{c_0} f_{\calG}(g) g^{p-1}dg \times \nu_{\calZ}(dz), 
\end{equation}
where $c_0=\int_0^\infty f_{\calG}(g) g^{p-1}dg$ 
and $\nu_{\calZ}$ is a probability measure on $\calZ=\{x : g(x)=1 \}.$ 

Let $dz'$ denote the volume element of the unit sphere 
$\bbS^{p-1} \subset \bbR^p.$ 
Then $c_0$ can also be written as 
$c_0= 1/\int_{\bbS^{p-1}}g(z')^{-p} dz',$ and the distribution of 
the direction $z'=x/\Vert x \Vert$ is given by 
$
   c_0\, g(z')^{-p} dz'. 
$

Under the additional assumption that $g(x)$ is 
piecewise of class $C^1,$ we can write the 
$\nu_{\calZ}(dz)$ in (\ref{eq:joint-star}) as 
$
  \nu_{\calZ}(dz) = c_0 \, \langle z, n_z \rangle \,dz,
$
where $n_z$ is the outward unit normal vector of
$\calZ$ and $dz$ on the right-hand side 
is the volume element of $\calZ.$  
\end{theorem}

In addition, Corollary \ref{co:given-density-on-Z'} 
together with Remark \ref{rem:given-density-on-Z'} 
yields the following result 
for the case of star-shaped distributions.

\begin{co} \quad
Suppose we are given an arbitrary distribution on $\bbS^{p-1}$ 
which has almost everywhere positive density $f(z')$ 
with respect to the volume element $dz'$ on $\bbS^{p-1}.$ 
Then we can realize this distribution $f(z')dz'$ as 
the distribution of the direction $z'=x/\Vert x \Vert$ of 
$x \in \bbR^p-\{ {\bf 0} \}$ which is distributed according to 
a star-shaped distribution.
\end{co}
\Proof
Since $g(z')^{-1} \in \Delta^{-1}(\{ f(z') \})=\{ f(z')^{1/p}\},$ 
we may take $\calZ = 
\left\{ f(z')^{1/p} z' : z' \in \bbS^{p-1} \right\}$ and 
$g(x) = \Vert x \Vert g (x/\Vert x \Vert) 
= \Vert x \Vert f(x/\Vert x \Vert)^{-1/p}.$
\QED

\section{Applications to random matrices}
\label{sec:random-matrix}

In this section we consider cross-sectionally contoured 
distributions of random matrices.  
For illustrative purposes, we consider a generalization of
matrix beta distribution by taking actions of the 
triangular group and the general linear group.  
These groups are not commutative.
Furthermore, the action of the general linear group is not free.  
Therefore, the results of Sections \ref{sec:cross-section} 
and \ref{sec:decomposable-distribution} 
can be fully illustrated by this example.  
Other examples of decomposable distributions of random matrices are 
given in \cite{Takemura-Kuriki96} and \cite{Kamiya-Takemura}.
See also \cite{Cadet} for a generalization of elliptically contoured
distribution to random matrices.

Let $W_1=(w_{1,ij})$ and $W_2=(w_{2,ij})$ be two $p\times p$ 
positive definite matrices. 
The sample space $\calX$ is $\{W=(W_1, W_2): W_1,W_2 \in PD(p)\}$ 
(Section 5.1) or essentially this set but with some exceptional 
null subset removed (Section 5.2).

As a dominating measure on $\cal X,$ we consider 
\begin{equation}
\label{eq:w1w2-measure}
\lambda(dW)=(\det W_1)^{a-\frac{p+1}{2}}
(\det W_2)^{b-\frac{p+1}{2}} dW_1 dW_2, 
\end{equation}
where $a,b > (p-1)/2$ and 
$dW_1 = \prod_{1 \le i \le j \le p} dw_{1,ij}, \ 
 dW_2 = \prod_{1 \le i \le j \le p} dw_{2,ij}.$

\subsection{Action of the triangular group}

First we consider the action of the lower triangular group.  
Let $LT(p)$ denote the group consisting  
of $p \times p$ lower triangular 
matrices with positive diagonal elements. 
Then $\calG=LT(p)$ acts on 
\[
\calX = \{(W_1, W_2): W_1,W_2 \in PD(p)\}
\] 
by
\[
(T, \ (W_1,W_2)) \mapsto (TW_1T\trans, \ TW_2 T\trans), \ \ T \in LT(p).
\]
This action is free and any cross section under this action 
is global.

It is interesting to note that there are two common cross sections 
used in the literature.  
Let $TT\trans =W_1 + W_2$ be the Cholesky 
decomposition of $W_1+W_2.$  
Then $T$ itself is an equivariant function and
$U = T^{-1} W_1 (T^{-1})^\trans$ is the associated invariant function.  
If $W_1$ and $W_2$ are independent Wishart matrices, then
$U$ has the matrix beta distribution.  
On the other hand, let $TT\trans=W_2$ be 
the Cholesky decomposition of $W_2.$
Then the invariant $F=T^{-1} W_1 (T^{-1})\trans$ has the 
matrix F distribution (\cite{Dawid}, %
Chapter 5 of \cite{Farrell}).

Here we prefer to consider the Cholesky decomposition of $W_1+W_2$
and use the following beta-type cross section:
\[
\calZ' = \left\{ (U, \ I_p-U): %
0 < U < I_p \right\}
\subset PD(p)\times PD(p),
\]
where $I_p$ denotes the $p \times p$ identity matrix 
and $A<B$ means $B-A \in PD(p)$ for $p\times p$ symmetric 
$A$ and $B.$ 
The orbital decomposition of $W=(W_1, W_2)$ with respect to $\calZ'$ 
is written as  
\begin{equation}
\label{eq:Z'triangular}
 (W_1, W_2) = \left(TUT\trans, \ T(I_p-U)T\trans\right), 
\qquad T=T(W),\ U=U(W).
\end{equation}

Next we move on to a general cross section. 
By using Remark \ref{rem:variety(free)} in the opposite direction, 
we obtain a general cross section $\calZ:$  
\[
\calZ = \left\{ z_U :  
0 < U < I_p \right\} 
\]
with 
\[
z_U = \left(S(U)US(U)\trans, \ S(U)(I_p - U) S(U)\trans \right), 
\]
where $S(U)=\left(s_{ij}(U)\right)$ is a function 
from $\left\{ U : %
0 < U < I_p \right\}$ 
to $LT(p).$  
Then the associated equivariant function is 
\[
g(W) = T(W)S(U(W))^{-1} 
\]
by Remark \ref{rem:transformation of equivariant part (free)}, 
and the invariant part is 
\[
z(W)=z_{U(W)}
=\left(S(U(W))U(W)S(U(W))\trans, \ 
S(U(W))(I_p - U(W)) S(U(W))\trans \right).
\]

Using a density of the form
\begin{equation*}
f(W) %
= f_{\calG}( g(W) )
\end{equation*}
with respect to $\lambda(dW)$ in (\ref{eq:w1w2-measure}), 
we obtain a cross-sectionally contoured distribution 
with respect to $\calZ.$ 
Now the application of 
Theorems \ref{th:distributions of y and z} 
and \ref{th:distribution of z^prime}
gives the following results about the distributions of 
$G = g(W), \ Z = z(W)$ 
and $U=U(W).$ 
Note that $U(W)$ is in one-to-one correspondence with 
the invariant part $z'(W)=\left(U(W), \ I_p-U(W) \right)$ with 
respect to $\calZ'.$

\begin{theorem} \quad
\label{Wishart-2-sample}
Suppose that the distribution of $W=(W_1,W_2)$ is given as 
\[
f_{\calG}(g(W)) (\det W_1)^{a-\frac{p+1}{2}} 
(\det W_2)^{b -\frac{p+1}{2}} dW_1 dW_2
\]
with some $f_{\calG}: LT(p) \rightarrow \bbR.$
Then $G=(g_{ij})=g(W)$ and $Z=z(W)$ are independent, 
and their joint distribution is given by 
\begin{equation*}
  \label{density-for-Wishart-2-sample-cross-sectionally-contoured}
\frac{1}{c_0} \,f_{\calG}(G) \prod_{i=1}^p g_{ii}^{2(a+b)-i} dG \times  
\nu_{\calZ}(dZ),
\end{equation*}
where
$
c_0 = \int_{LT(p)} f_{\calG}(G) \prod_{i=1}^p g_{ii}^{2(a+b)-i} dG 
$ 
and $\nu_Z$ is a probability measure on $\calZ.$ 
Furthermore, $U=(u_{ij})=U(W)$ is independent of $G=g(W),$ and its 
distribution is given by 
\[
2^p c_0 \prod_{i=1}^p s_{ii}(U)^{2(a+b)+p-2i+1}
(\det U)^{a-\frac{p+1}{2}} (\det (I_p-U))^{b-\frac{p+1}{2}} dU,
\]
where $dU=\prod_{1 \le i \le j \le p} du_{ij}.$
\end{theorem}
\Proof 
Remember the following well-known facts:
(a) The multiplier of relatively invariant measure $\lambda(dW)$ in 
(\ref{eq:w1w2-measure}) is 
$\chi(T)=(\det T)^{2(a+b)}$ 
(\cite{Wijsman90}, (9.1.4)), 
so $\chi(T)= \prod_{i=1}^p t_{ii}^{2(a+b)}$  
for $T=(t_{ij})\in LT(p);$
(b) For $LT(p),$ the left Haar measure is a multiple of  
$
\prod_{i=1}^p t_{ii}^{-i} dT 
$ 
(\cite{Wijsman90}, (7.7.2)) 
and the right-hand modulus is 
$
\Delta(T)^{LT(p)}= \prod_{i=1}^p t_{ii}^{p-2i+1}
$
(\cite{Wijsman90}, (7.7.6));
(c) With respect to the standard cross section $\calZ',$ we have the 
factorization 
$dW_1 dW_2 %
            = 2^p \prod_{i=1}^{p} t_{ii}^{2p+2-i} dT dU,
 \ T=(t_{ij}),$ for (\ref{eq:Z'triangular}) 
(\cite{Farrell}, (10.3.5)).
With the help of these facts,  
the theorem follows immediately from 
Theorems \ref{th:distributions of y and z} 
and \ref{th:distribution of z^prime}.
\QED

\subsection{Action of the general linear group}
\label{subsec:two-wishart-gl}

Consider the action of the general linear group $\calG=GL(p)$ 
consisting of all $p\times p$ nonsingular matrices.  
In this case, for there to exist a global cross section, 
we restrict the sample space as
\begin{eqnarray}
&& \calX= \{ (W_1, W_2) \in PD(p) \times PD(p): \nonumber\\
&& \qquad \qquad 
   \mbox{the $p$ roots of}\ \det(W_1- l(W_1+W_2))=0 \ 
   \mbox{are all distinct} \}.
\label{eq:wishart sample space}
\end{eqnarray}
If there are multiple roots in (\ref{eq:wishart sample space}), 
there are more than one orbit type (Appendix \ref{app:orbit-types}).
As in the case of $LT(p),$ the action of $GL(p)$ is 
$(B, \ (W_1, W_2)) \mapsto (BW_1B\trans, \ BW_2B\trans), \ 
B \in GL(p).$

As a standard global cross section,  
we can take
\begin{equation*}
\label{eq:wishart cross section}
  \calZ'= \left\{ (L, I_p- L): 
  L = \diag (l_1,\ldots,l_p),
  \ 1 > l_1 > \cdots > l_p > 0 \right\}.
\end{equation*}
For this $\calZ',$ the common isotropy subgroup is 
\begin{equation}
\label{eq:G0'}
  \calG_0' = \left\{ \diag (\epsilon_1,\ldots,\epsilon_p):  
  \epsilon_1 = \pm 1,\ldots,\epsilon_p = \pm 1 \right\},
\end{equation}
and the normalizer of $\calG_0'$ is given as
\begin{eqnarray*}
\calN' = \{ P \in GL(p) &:&
P \text{ has exactly one nonzero element} \\
&{}& \text{in each row and in each column} \}, 
\end{eqnarray*}
which is the group generated by permutation matrices and 
nonsingular diagonal matrices.
The orbital decomposition of $W=(W_1,W_2)$ with respect to $\calZ'$ 
can be written as
\begin{eqnarray}
\label{eq:BLB'}
(W_1,W_2) &=& (B L B\trans, \ B(I_p-L)B\trans) \\
          &=& (B(W) L(W) B(W)\trans, \ 
               B(W)(I_p-L(W))B(W)\trans). \nonumber
\end{eqnarray}
In this representation, $L(W)=\diag(l_1(W),\ldots, l_p(W))$ 
is uniquely determined by $W,$ but $B(W)$ is unique only up to 
the sign of each column of $B(W).$ 
We can use an arbitrary selection $B(W),$ e.g., 
the selection $B(W)$ such that $B(W)\in GL(p)/2^p,$ where 
$GL(p)/2^p$ denotes the set of $p\times p$ nonsingular matrices 
whose first nonzero element in each column is positive.
Note that it seems more convenient here to work with a 
selection $B(W)$ rather than the cosets of ${\calG}_0'$,  
although they are equivalent. 

Now we turn to a general global cross section. 
Using Theorem \ref{th:variety} in the opposite direction, 
we find that 
a general global cross section $\calZ$ %
is of the form 

\begin{equation*}
\calZ = \left\{ z_{L}: 
L = \diag (l_1,\ldots,l_p),
  \ 1 > l_1 > \cdots > l_p > 0 \right\}, 
\end{equation*}
where 
\begin{equation}
\label{eq:arbitrary cross section for two-sample Wishart}
z_{L}=
\left(B_0 P(L) L P(L)\trans B_0\trans, \ 
B_0 P(L)(I_p-L) P(L)\trans B_0\trans \right) 
\end{equation}
with $B_0 \in GL(p)$ and $P(L) \in \calN'.$ 
Without loss of generality, we take $B_0=I_p$ in
(\ref{eq:arbitrary cross section for two-sample Wishart}) 
so that the isotropy subgroup for $\calZ$ is also the 
$\calG_0'$ in (\ref{eq:G0'}).
By Proposition \ref{prop:transformation of equivariant part}, 
a selection of the equivariant part with respect to 
$\calZ$ is given by 
\begin{equation}
\label{eq:BP(L)^-1}
B(W)P(L(W))^{-1},
\end{equation}
but here we take the selection $g(W)$ which is given by changing the sign 
of each column of (\ref{eq:BP(L)^-1}) if necessary so that 
$g(W)\in GL(p)/2^p.$ 
The invariant part, on the other hand, is 
\[
z(W)
= z_{L(W)} 
= \left(P(L(W)) L(W) P(L(W))\trans, \ 
P(L(W))(I_p-L(W)) P(L(W))\trans \right).
\]

Consider a density of the form
\begin{equation*}
\label{eq:example of nonstandard density for wishart}
f(W) %
= t( g(W) ), 
\end{equation*}
where $t: GL(p) \rightarrow \bbR$ 
satisfies $t(BB_1)=t(B), \ B \in GL(p), \ B_1 \in \calG_0'.$
Then $t( g(W) ) \lambda(dW)$ is  
a cross-sectionally contoured distribution with respect to $\calZ.$ 
Applying Theorems \ref{th:distributions of y and z} 
and \ref{th:distribution of z^prime}, we obtain 
the following results about the distributions of 
$G = g(W), \ Z = z(W)$ 
and $%
(l_1,\ldots,l_p)=(l_1(W),\ldots,l_p(W)).$ 
Notice that $(l_1(W), \ldots, l_p(W))$ is in one-to-one correspondence 
with $z'(W)=\left(L(W), \ I_p-L(W) \right).$

\begin{theorem} \quad
\label{Wishart-2-sample-GL}
Suppose that the distribution of $W=(W_1,W_2)$ is %
\[
t(g(W))(\det W_1)^{a-\frac{p+1}{2}}
(\det W_2)^{b-\frac{p+1}{2}}dW_1 dW_2, 
\]
where $t: GL(p) \rightarrow \bbR$ is a real-valued function such that 
$t(B), \ B \in GL(p),$ does not depend on the sign 
of each column of $B.$ 
Then $G = g(W)$ and $Z = z(W)$ are independent, 
and their joint distribution is given by 
\begin{equation*}
\label{density-for-Wishart-2-sample-cross-sectionally-contoured-GL}
\frac{1}{c_0} %
t(G) (\det G)^{2(a+b)-p} %
dG %
\times \nu_{\calZ}(dZ),
\end{equation*}
where 
$c_0 = \int_{GL(p)/2^p}t(G)(\det G)^{2(a+b)-p}dG
= 2^{-p}\int_{GL(p)}t(G)(\det G)^{2(a+b)-p}dG$
and $\nu_{\calZ}$ is a probability measure on ${\cal Z}.$
Furthermore, $(l_1,\ldots, l_p)=(l_1(W),\ldots, l_p(W))$ 
is independent of $G=g(W),$ and its distribution is given by 
\[
2^p c_0 \, 
(\det P(\diag(l_1,\ldots,l_p)))^{2(a+b)} \, 
\prod_{i=1}^p l_i^{a-\frac{p+1}{2}} \, 
\prod_{i=1}^p (1-l_i)^{b-\frac{p+1}{2}} \, 
\prod_{i<j} (l_i - l_j) \,
dl_1 \cdots dl_p.
\]
\end{theorem}
\Proof 
This theorem 
is a direct consequence of 
Theorems \ref{th:distributions of y and z} 
and \ref{th:distribution of z^prime}. 
We only have to recall the following easy or well-known facts:
(a) The (left) Haar measure $\mu_{GL(p)}$ on $GL(p)$ %
is a multiple of  
$(\det B)^{-p}dB;$ 
(b) $GL(p)$ is unimodular, 
so $\Delta(B)=\chi(B)=(\det B)^{2(a+b)}, \ B \in GL(p);$
(c) In terms of the standard global cross section $\calZ',$ 
we have the 
factorization $dW_1 dW_2 = 2^p (\det B)^{p+2}dB 
\prod_{i<j}(l_i-l_j) dl_1 \cdots dl_p$ 
for (\ref{eq:BLB'})
(\cite{Anderson2003}, Theorem 13.2.1). 
\QED

\appendix

\

\begin{flushleft}
{\bf \Large Appendix}
\end{flushleft}

\section{Variety of global cross sections}
\label{app:variety}

In this Appendix, we discuss the construction of general global cross
sections from a given global cross section and characterize the class
of all global cross sections in terms of the normalizer of the common
isotropy subgroup.  
The proof of Theorem \ref{th:variety} is provided in this Appendix 
in particular, but a thorough investigation 
into the variety of global cross sections, 
including the uniqueness of $n_z$ in (\ref{eq:g0nz}), 
is also conducted here.
This material was partly discussed in
\cite{Kamiya-Takemura}, but here we give a complete characterization.

\subsection{Action of a factor group on each $\calG$\kmms-orbit}

We begin by confining our discussion to the action of $\calG$ on each
$\calG$\kmms-orbit 
$
\tilde{ \calX } = \calG x_0%
$
with $x_0 \in \calX.$
For an arbitrary point $x \in \tilde{\calX},$ let 
\[
\calN_x := \{ g \in \calG: g\calG_x g^{-1} = \calG_x \}
\] 
be the normalizer of $\calG_x= \{ g \in \calG: gx=x \}$ in $\calG.$
Then, 
since $\calG_{gx}=g\calG_xg^{-1}$ for 
$g \in \calG$ and $x \in \tilde{\calX},$ we have 
that the normalizers satisfy  (\cite{Kawakubo}, p.33)
\[
\calN_{gx}=g\calN_xg^{-1}, 
\ \ \ g \in \calG, \ x \in \tilde{\calX}.
\]

Now, consider the factor group 
\[
\calM_x := \calN_x/\calG_x = \{ n\calG_x: n \in \calN_x \}
\]
for each $x \in \tilde{\calX}.$
Then we have the following proposition: 

\begin{prop} \quad
\label{prop:isomorphic}
All factor groups $\calM_x, \ x \in \tilde{\calX},$ 
are isomorphic to one another.
\end{prop}

\Proof
For given $\calM_x$ and $\calM_{x'}, \ x, x' \in \tilde{ \calX },$
take an element $g \in \calG$ such that $x' = gx.$ 
Then we have $\calG_{x'}=g \calG_x g^{-1},$
so Lemma 1.51 of \cite{Kawakubo} implies that 
the mapping 
\[
\tau_{x', x}: \calM_x \to \calM_{x'}
\] 
defined as
$m = n \calG_x \mapsto (gng^{-1})(g \calG_x g^{-1})
= gmg^{-1}, \ n \in \calN_x,$ 
serves as an isomorphism.
\QED

Fix $z \in \tilde{ \calX }$ as a reference point, and write
\[
\calG_0 = \calG_z, \ \ \calN = \calN_z 
\ \ {\rm and} \ \ \calM = \calM_z.
\]
Now we define an action of $\calM$ on
$\tilde{ \calX }$ as follows:

Let $\tau_x, \ x \in \tilde{ \calX },$ be an arbitrary selection
of $\tau_{x, z}$: 
\[
\tau_x( m ) \in \tau_{x, z}( m ) \subset \calN_x, 
\ \ \ m \in \calM.
\]
Using this $\tau_x,$ 
we define $xm, \ x \in \tilde{\calX}=\calG z, \ m \in \calM,$ 
as 
\[
xm := \tau_x(m)x.
\]
If we write $x=gz, \ g \in \calG$ and 
$m = n\calG_0, \ n \in \calN,$ 
we can express $xm$ as
\begin{equation}
\label{xm}
xm = \tau_x(m)x = (gng^{-1})(gz)=gnz.
\end{equation}

We can confirm that this is well-defined:
For $m = n\calG_0 = ng_0\calG_0, \ g_0 \in \calG_0,$ and
$x = gz = g g_0' z, \ g_0' \in \calG_0,$ we have
$xm = (gg_0')(ng_0)z = g g_0' n z = gng_0''z = gnz$ since $g_0'n = ng_0''$
for some $g_0'' \in \calG_0.$

Moreover, $xm, \ x \in \tilde{\calX}, \ m \in \calM,$ has %
the following property.

\begin{lm} \quad
\label{lm:xm1m2}
For any $x \in \tilde{\calX}$ and $m_1, m_2 \in \calM,$ we have
\[
(xm_1)m_2 = x(m_1m_2).
\]
\end{lm}

\Proof
For $x = gz \in \tilde{ \calX }, \ g \in \calG,$ and 
$m_i = n_i \calG_0 \in \calM, \ n_i \in \calN, \ i=1,2,$ we have
$(xm_1)m_2 = (gn_1z)m_2 = (gn_1)n_2z = g(n_1n_2)z = x(m_1m_2)$
since $m_1m_2 = (n_1n_2)\calG_0.$
\QED

So we have a right action of $\calM$ on $\tilde{\calX}$:
\begin{prop} \quad
\label{prop:right}
The mapping 
\[
(x, m) \mapsto xm, \ \ \ x \in \tilde{\calX}, \ m \in \calM,
\]
is a right action of $\calM$ on $\tilde{ \calX }.$
Moreover, this action is free.
\end{prop}

\Proof
It is easy to see that $e_{ \calM }$   
(the identity element of $\calM$\kmms) satisfies 
$x e_{ \calM }=x$ for all $x \in \tilde{\calX}.$
This fact together with lemma \ref{lm:xm1m2} implies 
that the mapping $(x, m) \mapsto xm$ is a right action.

Next we show that this action is free.
Suppose $xm = x$ for $m = n\calG_0, \ n \in \calN$ and  $x = gz, \ g \in \calG.$
Then, $gnz = gz$ or $nz = z.$ So $n \in \calG_0,$ and
$m = n\calG_0 = \calG_0 = e_{ \calM }.$
\QED

As the following proposition shows, the $\calM$\kmms-orbits 
$x\calM, \ x \in \tilde{\calX},$ 
can be characterized in terms of the isotropy subgroups 
$\calG_x, \ x \in \tilde{\calX},$ 
under the action of $\calG$ on $\tilde{\calX}.$

\begin{prop} \quad
\label{prop:M-orbit}
For $x, x' \in \tilde{ \calX },$ we have that
\[
x\calM =x'\calM
\ \ \ {\rm if \ and \ only \ if} \ \ \ 
\calG_x = \calG_{x'}.
\]
\end{prop}

\Proof
Suppose 
$x\calM =x'\calM.$
Then, writing $x' = xm, \ x =gz, \ g \in \calG, \
m = n\calG_0, \ n \in \calN,$ we can calculate
$\calG_{x'} = \calG_{xm} = \calG_{gnz} = gn \calG_0 n^{-1}g^{-1}
= g \calG_z g^{-1} = \calG_{gz} = \calG_x.$

Conversely, suppose $\calG_x = \calG_{x'}.$
Then, writing $x' = g'x, \ g' \in \calG, \ x = gz, \ g \in \calG,$
we have $\calG_{x'} = \calG_{g'x} = g' \calG_x g'^{-1} = \calG_x,$
and thus $g' \in \calN_x = g \calN g^{-1}.$
Accordingly, we can write $g' = gng^{-1}, \ n \in \calN,$
and hence $x' = g'x = (gng^{-1})(gz) = gnz = xm$ for
$m = n\calG_0 \in \calM.$
Therefore, 
$x\calM = x'\calM.$
\QED

So the isotropy subgroups 
$\calG_x, \ x \in \tilde{\calX},$ are 
constant on each $\calM$\kmms-orbit and different on different 
$\calM$\kmms-orbits.
Therefore, 
the $\calM$\kmms-orbits can be labeled by the 
isotropy subgroups $\calG_x, \ x \in \tilde{\calX},$ which 
do not depend on the choice of the reference point $z.$

Now we have two groups $\calG$ and $\calM$ acting 
on $\tilde{\calX}.$
These two actions commute with each other:

\begin{prop} \quad
\label{prop:gxm}
The actions of $\calG$ and $\calM$ on $\tilde{\calX}$ commute:
\begin{equation}
\label{gxm}
g(xm) = (gx)m, \ \ \ g \in \calG, \ m \in \calM, \ x \in \tilde{ \calX }.
\end{equation}
\end{prop}

\Proof
Writing $x = g'z, \ g' \in \calG,$ and $m = n\calG_0, \ n \in \calN,$
we can deduce
$g(xm) = g(g'nz) = (gg')nz = (gg'z)m = (gx)m.$
\QED

Thus, we can say 
that $\tilde{\calX}$ is a $\calG$\kmms-\kmms$\calM$ bispace.
We will write $g(xm)$ and $(gx)m$ as $gxm$ without 
ambiguity.

As we saw in Proposition \ref{prop:M-orbit}, the $\calM$\kmms-orbits 
$x\calM, \ x \in \tilde{\calX},$ can be labeled without 
referring to the reference point $z,$ but the 
$\calM=\calM_z$ itself does depend on $z.$
By using the commutativity of the actions of $\calG$ and 
$\calM$ in Proposition \ref{prop:gxm}, 
and considering proportional translations by the action of $\calG,$ 
we can identify the elements of $\calM$ in terms of %
relative positions of two points of $\tilde{\calX}$ 
and thereby get rid of $z$ as follows: 

Since the action of $\calM$ on $\tilde{\calX}$ is free 
by Proposition \ref{prop:right}, we know that $\calM$ can be 
identified with an $\calM$\kmms-orbit: 
\[
\calM \leftrightarrow x\calM, \ \ \ x \in \tilde{\calX}.
\]
So we can see $\calM$ as $z\calM$ in particular, and hence  
as $\{ z \} \times z\calM$ 
by $m \leftrightarrow (z, zm):$ 
\[
\calM \leftrightarrow \{ z \} \times z\calM, \ \ \ 
m \leftrightarrow (z, zm), \ m \in \calM.
\]
Indicate by $\Pi$ 
the set of ordered pairs of points of $\tilde{\calX}$ on the same $\calM$\kmms-orbits:
\[
\Pi = \{ (x_1, x_2)\in \tilde{\calX}\times \tilde{\calX}: 
x_1\calM = x_2\calM \},
\]
and define an equivalence relation $\sim_{\Pi}$ among the 
elements of $\Pi$ as follows: 
\[
(x_1, x_2) \sim_{\Pi} (x_1', x_2')
 \ \Leftrightarrow \
g(x_1, x_2) = (x_1', x_2') 
{\rm \ \ for \ \ some \ \ }
g \in \calG,
\] 
where $g(x_1, x_2)$ is the proportional translate 
of $(x_1, x_2)$ by $g:$ 
\[
g(x_1, x_2) := (gx_1, gx_2), 
\]
i.e., the diagonal action. %
We denote the equivalence class under $\sim_{\Pi}$ by $[ \ \cdot \ ]_{\Pi}.$
Then, we can think of $\Pi/\kmms\sim_{\Pi}$ as a group with the following
product:
\begin{equation}
\label{operation}
[(x_1, x_2)]_{\Pi} \cdot [(x_3, x_4)]_{\Pi} = [(x_1, gx_4)]_{\Pi},
\end{equation}
where $g$ is an arbitrary element of $\calG$ satisfying $x_2 = gx_3.$
That is, $[(x_1, x_2)]_{\Pi} \cdot [(x_2, gx_4)]_{\Pi} = [(x_1, gx_4)]_{\Pi}.$
We can check that operation (\ref{operation}) is well-defined. 

\begin{prop} \quad
\label{prop:M=Pi/sim}
The factor group $\calM$ is isomorphic to the group $\Pi /\kmms\sim_{\Pi}.$
\end{prop}

\Proof
Consider the following two mappings:
\begin{eqnarray}
\label{first}
\calM \ni m &\mapsto& [(z, zm)]_{ \Pi } \in \Pi/\kmms\sim_{\Pi}, \\
\label{second}
\Pi/\kmms\sim_{\Pi} \ni [(x_1, x_2)]_{ \Pi } &\mapsto& {\rm the \ unique \ } m \in \calM 
{\rm \ such \ that \ } x_2 = x_1m.
\end{eqnarray}
It can be verified that (\ref{first}) and (\ref{second}) are 
the inverse mappings of each other. 
Moreover, we can show that (\ref{first}) is a homomorphism 
in the following way: 
For $m_1, m_2 \in \calM,$ we have 
$[(z, z(m_1m_2))]_{\Pi}=[(z, (n_1n_2)z)]_{\Pi}
=[(z, n_1(z m_2)]_{\Pi}
=[(z, z m_1)]_{\Pi} \cdot [(z, z m_2)]_{\Pi}$ with $n_1 \in m_1$ and 
$n_2 \in m_2.$ 
\QED

Note that thanks to Proposition \ref{prop:M-orbit}, 
the group $\Pi/\kmms\sim_{ \Pi }$ does not 
depend on the reference point $z.$
Proposition \ref{prop:M=Pi/sim} implies that an element of 
$\calM$ can be specified by an ordered pair $(x_1, x_2)$ of points of 
$\tilde{\calX}$ 
having the same isotropy subgroup $\calG_{x_1}=\calG_{x_2}$ 
if we identify all proportional translates $g(x_1, x_2), \ g \in \calG.$

\subsection{Global cross sections on the whole sample space}

Let us get back to the action of $\calG$ on the whole of $\calX.$
Throughout this subsection, we assume that there exists 
a global cross section. 
We agree that a global cross section always refers to 
the one under the action of $\calG$ on $\calX.$

For an arbitrary global cross section $\calZ',$ we write
\[
\calG_{0, \calZ'} = \calG_{z'}, \ \ \calN_{\calZ'} = \calN_{z'}
\ \ {\rm and} \ \ 
\calM_{\calZ'} = \calM_{z'}
\] 
with $z' \in \calZ'.$

First we note that 
the difference between two proportional global cross sections 
$\calZ'$ and $g\calZ'$ is not essential since 
they induce the same family of proportional global cross sections.
So we introduce the equivalence relation $\sim_{{\rm gcs}}$ 
among the global cross sections by proportionality:
\[
  \calZ_1 \sim_{{\rm gcs}} \calZ_2 \ \Leftrightarrow \
  \calZ_1 = g\calZ_2 
  {\rm \ \ for \ \ some \ \ }
  g \in \calG.
\]
The equivalence class under $\sim_{{\rm gcs}}$ 
is indicated by $[ \ \cdot \ ]_{{\rm gcs}}.$

Fix $\calZ$ as a reference global cross section, and put
\[
\calG_0 = \calG_{0, \calZ}, \ \ \calN = \calN_{ \calZ }
\ \ {\rm and} \ \ 
\calM = \calM_{ \calZ }.
\]
For this $\calZ,$ let 
$\iota = \iota_{\calZ},$
i.e., 
the natural one-to-one correspondence between $\calX/\calG$ 
and $\calZ.$
Based on the reference global cross section $\calZ,$ 
we can generate a global cross section from another global cross 
section in the following way.

Let $\calZ'$ be an arbitrary global cross section,
and $M$ a mapping from $\calX / \calG$ to $\calM.$
Consider the following subset of $\calX:$  
\begin{equation}
\label{Z'm}
\calZ' M := \{ z' m_{z'}: z' \in \calZ' \}, 
\end{equation}
where 
$m_{x}:= M( \calG x )$ for $x \in \calX.$
Note that 
\[
m_{gx}=m_{x}, \ \ \ g \in \calG, \ x \in \calX.
\] 
In (\ref{Z'm}), $z'm_{z'}$ is defined as in (\ref{xm}) under
the action of $\calM$ on $\tilde{\calX} = \calG z'$ with
$z = \iota(\calG z') \in \calG z' \cap \calZ$  
as the reference point of $\calG z'.$
If we write $z'=g_z z, \ g_z \in \calG,$ 
and use an arbitrary 
$n_z \in m_z = m_{z'} \in \calM$ for each $z'\in \calZ',$ 
we obtain a more direct expression of definition (\ref{Z'm}): 
\[
\calZ' M = \{ g_z n_z z: z \in \calZ \}.
\]
Now we have that $\calZ'M$ is also a global cross section:

\begin{theorem} \quad
\label{th:Z'm}
Let $\calZ'$ be an arbitrary global cross section, 
and let $M$ be an arbitrary mapping from $\calX / \calG$ to $\calM.$
Then, $\calZ' M$ is a global cross section.
\end{theorem}

\Proof
It is clear that $\calZ'M =\{ z'm_{z'}: z' \in \calZ' \}$ 
is a cross section, 
since $z'$ and $z'm_{z'}\in z'\calM \subset \calG z'$ 
are on the same $\calG$\kmms-orbit $\calG z'.$
Further, the cross section $\calZ'M$ is global because 
for $z'=g_z z$ and $n_z \in m_{z'},$ we have 
$\calG_{g_z n_z z}
=g_z n_z \calG_0 n_z^{-1}g_z^{-1}
=g_z \calG_z g_z^{-1}
=\calG_{g_z z}
=\calG_{z'}
=\calG_{0, \calZ'},$ 
common for all $z\in \calZ.$ 
\QED

Theorem \ref{th:Z'm} implies, in particular, that 
for arbitrary $g \in \calG$ and $M: \calX / \calG \to \calM,$   
the subset $(g \calZ)M \subset \calX$ is a global cross section.
We will show below that the converse is true as well.
That is, an arbitrary global cross section $\calZ'$ must be 
of this form: 
\[
\calZ'=(g \calZ)M, \ \ \ g \in \calG, \ \ M: \calX / \calG \to \calM.
\]
Moreover, we want to study the uniqueness of $M$ in such an 
expression.
To make the arguments succinct, 
we will introduce an action on the set of 
equivalence classes of global cross sections as follows.

Using $\calZ'M$ in (\ref{Z'm}), 
we define $[\calZ']_{{\rm gcs}}M$ as 
\[
[\calZ']_{{\rm gcs}}M := [\calZ' M]_{{\rm gcs}}
\]
for a global cross section $\calZ'$ and a mapping 
$M: \calX/\calG \to \calM.$
Note that $\calZ' M$ is a global cross section because of
Theorem \ref{th:Z'm}.
Thanks to the next lemma, this $[\calZ']_{{\rm gcs}}M$ 
is well-defined.

\begin{lm} \quad
\label{lm:gZ'm}
Let $\calZ'$ be an arbitrary global cross section.
Then, we have
\[
(g\calZ')M = g(\calZ'M)
\]
for any $g \in \calG$ and $M: \calX / \calG \to \calM.$
\end{lm}

\Proof
This is a direct consequence of Proposition \ref{prop:gxm}:
Taking $x = z' \in \calZ'$ and $m = m_{gz'} = m_{z'}$
in (\ref{gxm}), we obtain $(gz')m_{gz'} = g(z'm_{z'}).$
The result follows from (\ref{Z'm}).
\QED

By this lemma, 
we can write $(g\calZ')M$ and $g(\calZ'M)$ as $g\calZ'M.$

Now, the set of mappings $M: \calX/\calG \to \calM$ 
forms a group with the product defined pointwise.
The identity element is the constant map onto $e_{\calM},$
and the inverse elements are pointwise inverses in $\calM.$

Viewing the set of mappings $M: \calX/\calG \to \calM$ in this way,
we can regard 
\begin{equation}
\label{[Z']m}
([\calZ']_{{\rm gcs}}, M) \mapsto [\calZ']_{{\rm gcs}}M
\end{equation} 
as a right action of the group of mappings 
$M: \calX/\calG \to \calM$ on the set of equivalence classes
of global cross sections $[ \calZ' ]_{{\rm gcs}}.$
This can be seen from the following lemma.

\begin{lm} \quad
\label{lm:Z'm'm''}
For an arbitrary global cross section $\calZ'$ and
mappings $M_1, M_2: \calX/\calG \to \calM,$
we have 
\[
([\calZ']_{{\rm gcs}}M_1)M_2 = [\calZ']_{{\rm gcs}}(M_1 M_2).
\]
\end{lm}
\Proof
This follows essentially from Lemma \ref{lm:xm1m2}: 
Denoting $m_{i, x}:=M_i(\calG x)$, $x \in \calX$, for 
$i=1,2,$ we can write 
$\calZ' M_1=\{ z'm_{1, z'}: z' \in \calZ' \}$ and hence  
$(\calZ' M_1)M_2=\{ (z'm_{1, z'})m_{2, z'm_{1, z'}}: 
z' \in \calZ' \}.$ 
But since $m_{2, z'm_{1, z'}}=m_{2, z'},$ we obtain 
$(\calZ' M_1)M_2=\{ (z'm_{1, z'})m_{2, z'}: z' \in \calZ' \}
=\{ z'(m_{1, z'}m_{2, z'}): z' \in \calZ' \}=\calZ'(M_1M_2)$ 
by Lemma \ref{lm:xm1m2}.  
\QED

We can show that action (\ref{[Z']m}) is transitive as follows:
Let $\calZ'$ be an arbitrary global cross section.
Then, since $\calZ'$ and $\calZ$ are cross sections, we can write
$\calZ' = \{ g_z z: z \in \calZ \}$ for some $g_z$\kmms s in $\calG.$
Fix an arbitrary $z_0 \in \calZ$ and then we have
$\calG_{g_z z} = \calG_{g_{z_0}z_0} = g_z \calG_0 g_z^{-1} = g_{z_0} \calG_0 g_{z_0}^{-1}$
or $g_{z_0}^{-1}g_z \in \calN$ for all $z \in \calZ.$
Putting $g = g_{z_0},$ we can represent $g_z = g n_z$ for certain   
$n_z$\kmms s in $\calN,$ and thus obtain
\begin{equation}
\label{Z'}
\calZ' = \{ g n_z z: z \in \calZ \}.
\end{equation}
Now define $M: \calX/\calG \to \calM$ as 
$M(\calG z) 
= n_z \calG_0 \in \calM, \ z \in \calZ.$
With this $M,$ we can express (\ref{Z'}) as 
$\calZ' = g \calZ M$ and then
arrive at 
\begin{equation}
\label{Z'=ZM}
[ \calZ']_{{\rm gcs}} = [ \calZ ]_{{\rm gcs}}M.
\end{equation}
This proves that action (\ref{[Z']m}) is transitive.

On the other hand, (\ref{[Z']m}) is not a free action.
However, we have a certain kind of uniqueness 
of $M$ in (\ref{Z'=ZM}).
To state this fact about uniqueness,
we introduce an equivalence relation $\sim_m$ among mappings
$M: \calX/\calG \to \calM$ in the following way:
\[
M_1 \sim_m M_2 \ \Leftrightarrow \ 
M_1( \ \cdot \ ) \equiv \bar{m} M_2( \ \cdot \ ) 
{\rm \ \ for \ \ some \ \ } \bar{m} \in \calM.
\]
We denote the equivalence class under $\sim_m$ by
$[ \ \cdot \ ]_m.$
Now we verify the uniqueness of $[M]_m$ for 
$M$ in (\ref{Z'=ZM}):

Suppose $M': \calX/\calG \to \calM$ also satisfies
$[\calZ']_{{\rm gcs}} = [\calZ]_{{\rm gcs}}M'.$
Then we have $[\calZ M]_{{\rm gcs}} = [\calZ M']_{{\rm gcs}}$ 
and thus
\begin{equation}
\label{n'zz}
\forall z \in \calZ: n_z z = g n'_z z
\end{equation}
for some $g \in \calG,$ 
where $n_z \in M(\calG z)$ and 
$n_z' \in M'(\calG z).$
By considering the isotropy subgroup at $n_z z = g n'_z z,$
we can easily see that the above $g$ is in $\calN.$
It follows from (\ref{n'zz}) that $n_z \calG_0 = g n'_z \calG_0$ 
or
\[
M(\calG z)= \bar{m} M'(\calG z)
{\rm \ \ for \ \ all \ \ } 
z \in \calZ
\]
with $\bar{m} := g \calG_0 \in \calM,$
and this means $[M]_m = [M']_m.$

Thus we have proved the following theorem.

\begin{theorem} \quad
\label{th:Z'=Zm'}
Let $\calZ'$ be an arbitrary global cross section.
Then, there exists a mapping $M: \calX/\calG \to \calM$ such that
\[
[\calZ']_{{\rm gcs}} = [\calZ]_{{\rm gcs}}M.
\]
Here, $[M]_m$ is uniquely determined.
\end{theorem}

Statements in Theorem \ref{th:Z'=Zm'} can be 
translated into relations between two arbitrary 
global cross sections $\calZ_1$ and $\calZ_2.$ 
To do this,  we need to introduce some notation.

For arbitrary mappings $M: \calX/\calG \to \calM$ 
and global cross sections $\calZ',$   
we write 
\[
M^{\calZ'}:=M' M
\] 
with $M': \calX/\calG \to \calM$ 
such that  
$[\calZ']_{{\rm gcs}} = [\calZ]_{{\rm gcs}}M'.$
Note that because of the uniqueness part of 
Theorem \ref{th:Z'=Zm'}, 
$[M^{\calZ'}]_m$ does not depend on the choice of such $M'.$

\begin{co} \quad
\label{co:Z2=Z1m12}
Let $\calZ_1, \ \calZ_2$ be arbitrary global cross sections.
Then, there exists a mapping $M_{12}: \calX/\calG \to \calM$ 
such that
\[
[\calZ_2]_{{\rm gcs}} = [\calZ_1]_{{\rm gcs}}M_{12}.
\]
Here, $M_{12}$ is unique in the sense that for two such $M_{12}$ and
$M_{12}',$ we have 
\[
[(M_{12})^{\calZ_1}]_m=[(M'_{12})^{\calZ_1}]_m.
\]
\end{co}

\Proof
As was shown above, 
action (\ref{[Z']m}) is transitive. 
This implies the existence of $M_{12}.$

Now we show the uniqueness of $M_{12}.$
Suppose there exists another $M_{12}': \calX/\calG \to \calM$ 
such that
$[\calZ_2]_{{\rm gcs}} = [\calZ_1]_{{\rm gcs}}M_{12}'.$
Then we have $[\calZ_1]_{{\rm gcs}}M_{12} 
= [\calZ_1]_{{\rm gcs}}M_{12}'.$
But by virtue of Lemma \ref{lm:Z'm'm''}, this can be written as
$[\calZ]_{{\rm gcs}}(M_1M_{12}) = [\calZ]_{{\rm gcs}}(M_1M_{12}')$ 
with $M_1$ satisfying $[\calZ_1]_{{\rm gcs}} 
= [\calZ]_{{\rm gcs}}M_1.$
Now, the uniqueness part of Theorem \ref{th:Z'=Zm'} yields
$[M_1M_{12}]_m = [M_1M_{12}']_m.$
\QED

\begin{co} \quad
\label{co:Zj=Zimij}
Suppose two global cross sections $\calZ_1, \ \calZ_2$ are
related by
\[
[\calZ_2]_{{\rm gcs}} = [\calZ_1]_{{\rm gcs}}M_{12}, \ \ \  
[\calZ_1]_{{\rm gcs}} = [\calZ_2]_{{\rm gcs}}M_{21}
\]
for $M_{12}, M_{21}: \calX/\calG \to \calM.$
Then, these $M_{12}$ and $M_{21}$ are the inverse elements 
of each other in the sense that 
\[
[(M_{12})^{\calZ_1}]_m = [(M_{21}^{-1})^{\calZ_1}]_m, \ \ \ 
[(M_{21})^{\calZ_2}]_m = [(M_{12}^{-1})^{\calZ_2}]_m.
\]
\end{co}

\Proof
Since $[\calZ_2]_{{\rm gcs}} = [\calZ_1]_{{\rm gcs}}M_{12} 
= [\calZ_1]_{{\rm gcs}}M_{21}^{-1},$ 
we have directly from Corollary \ref{co:Z2=Z1m12} that 
$[(M_{12})^{\calZ_1}]_m = [(M_{21}^{-1})^{\calZ_1}]_m.$
The relation $[(M_{21})^{\calZ_2}]_m = [(M_{12}^{-1})^{\calZ_2}]_m$
is shown in a similar manner.
\QED

Now we have a characterization of global cross sections. 
Theorem \ref{th:Z'=Zm'} implies that an arbitrary global 
cross section $\calZ'$ must be written as 
$\calZ'=g\calZ M$ for some $g \in \calG$ and 
$M: \calX/\calG \to \calM.$
Together with the remark after 
Theorem \ref{th:Z'm}, this leads to the following 
characterization:
\begin{co} \quad
\label{co:characterization}
A subset $\calZ' \subset \calX$ is a global cross section 
if and only if it is of the form 
\[
\calZ'=g\calZ M, \ \ \ g \in \calG, \ \ M: \calX/\calG \to \calM.
\]
\end{co}

Therefore, $\calZ'$ is a global cross section if and only if 
it can be written as 
\begin{equation*}
\label{Z'=gnz}
\calZ'=\{ g n_z z: z \in \calZ \}
\end{equation*}
for some $g \in \calG$ and $n_z \in \calN, \ z \in \calZ.$ 

\section{Orbit types}
\label{app:orbit-types}

In this paper we have discussed properties of a global cross section,
assuming one exists.  However, a global cross section does not always
exist.  
In discussing star-shaped distributions in Section
\ref{sec:star-shaped}, we omitted the origin from the sample space to
guarantee the existence of a global cross section.  
Similarly, in Section
\ref{subsec:two-wishart-gl} we assumed the distinctness of the 
roots of the characteristic equation.  
In these examples, the excluded sets are of
measure zero and can be ignored.  
However, there are some cases where
``singular sets'' have positive measure and can not be ignored.  
An example of this case is given
by the orthogonal projection of a random matrix $U$ onto the cone 
of nonnegative definite matrices
(\cite{Kuriki}, 
\cite{Kuriki-Takemura00}).

A global cross section exists if and only if all the 
isotropy subgroups
${\cal G}_x, \ x \in \calX,$ are conjugate to one another.
One can confirm this easily by recalling 
property (\ref{eq:G_gx}).

Now define the equivalence relation 
$\sim_{\calX}$ in $\calX$ 
by the conjugacy of the isotropy subgroups:
\begin{equation*}
\label{sim_X}
x \sim_{\calX} x' \ \Leftrightarrow \
\calG_x = g \calG_{x'} g^{-1} {\rm \ \ for \ \ some \ \ }
g \in \calG.
\end{equation*}
Then, even when a global cross section does not exist for the action
of $\calG$ on the whole of $\calX,$ there does exist a global 
cross section if we restrict our attention to 
the action of $\calG$ on each equivalence class under $\sim_{\calX}.$  
These equivalence classes are called the {\it orbit types}. 
See Section 1.8 of \cite{Kawakubo} or Section 1.4 of \cite{Bredon}.
We assume that the number of orbit types is at most countable.   
It is known that if $\calG$ is compact, the number of 
orbit types is actually finite (see Section 4.1 of \cite{Bredon}).

Let $\{ {\cal X}_i: i \ge 1 \}$ be the partition of ${\cal X}$
into orbit types.
By restricting the action $({\cal G}, {\cal X})$ of ${\cal G}$ on
${\cal X},$ we obtain the action $({\cal G}, {\cal X}_i)$ of
${\cal G}$ on each ${\cal X}_i, \ i \ge 1.$
For each $i \ge 1,$ let ${\cal Z}_i$ be a global cross section
for $({\cal G}, {\cal X}_i),$
and denote by ${\cal G}_i$ the common isotropy subgroup at 
the points of ${\cal Z}_i.$
Then for each $i \ge 1,$ we have the orbital decomposition of
${\cal X}_i:$
\begin{eqnarray*}
& {\cal X}_i \leftrightarrow {\cal Y}_i \times {\cal Z}_i, & \\
& x_i \leftrightarrow (y_i, z_i), &
\ \ \ x_i = g_i z_i,
\ \ y_i = g_i {\cal G}_i \in {\cal Y}_i = {\cal G}/{\cal G}_i.
\end{eqnarray*}

Write ${\cal Y} = \bigcup_i {\cal Y}_i$ and
${\cal Z} = \bigcup_i {\cal Z}_i,$
and define the functions $y:{\cal X} \to {\cal Y}$ and
$z:{\cal X} \to {\cal Z}$ by
\begin{equation*}
\label{y_i,z_i}
y(x) = y_i(x), \ z(x) = z_i(x)  {\rm \ \ if \ } 
x \in {\cal X}_i, \ \ \ i \ge 1.
\end{equation*}
Note that ${\cal Z}$ is a cross section for $({\cal G}, {\cal X}).$

Concerning topological questions, we assume
1, 3 and 4 of Assumption \ref{assumption:regularity1} %
with $\calX, \ \calG_0, \ \calZ$ and 
$x \leftrightarrow (y, z)$ replaced by
$\calX_i, \ \calG_i, \ \calZ_i$ and 
$x_i \leftrightarrow (y_i, z_i),$ respectively.  
On $\calX_i$ we consider a dominating measure
$\lambda_i$ which is relatively invariant with multiplier $\chi_i.$
We note that $\calG$ is metrizable by 2 of 
Assumption \ref{assumption:regularity1} (\cite{Ash}, A5.16 Theorem).
We regard the elements of $\calY$ as subsets of $\calG.$
By endowing $\calY$ with the Hausdorff distance, we make $\calY$
a metric space.  
(For details, see Appendix A.3 of our technical report 
\cite{Kamiya-Takemura}.) 

Let $\lambda(dx)=\sum_i I_{{\cal X}_i}(x) \lambda_i(dx),$
where $I_{ \calX_i }$ is the indicator function of $\calX_i.$
Note that $\lambda$ is not necessarily a relatively invariant 
measure.
Now suppose that $x$ is distributed according to 
$
f_{\calY}(y(x))\lambda(dx)
$
for some $f_{\calY}:{\cal Y} \to \bbR.$
Here we assume
$\int_{{\cal X}_i} f_{\calY}(y(x)) \lambda(dx) > 0$
for each $i \ge 1.$
Under these conditions, it is easy to show that for each $i \ge 1,$
\begin{equation*}
\label{conditional independence}
P(y(x) \in A, \ z(x) \in B \mid x \in {\cal X}_i)
= P(y(x) \in A \mid x \in {\cal X}_i)
P(z(x) \in B \mid x \in {\cal X}_i)
\end{equation*}
for each measurable $A \subset {\cal Y}_i$ and $B \subset {\cal Z}_i.$
Therefore, $y(x)$ and $z(x)$ are conditionally
independent given $x \in \calX_i.$

\end{document}